# THE EULER SCHEME FOR LÉVY DRIVEN STOCHASTIC DIFFERENTIAL EQUATIONS: LIMIT THEOREMS

By Jean Jacod

*Université Pierre et Marie Curie*

We study the Euler scheme for a stochastic differential equation driven by a Lévy process $Y$. More precisely, we look at the asymptotic behavior of the normalized error process $u_n(X^n - X)$, where $X$ is the true solution and $X^n$ is its Euler approximation with stepsize $1/n$, and $u_n$ is an appropriate rate going to infinity: if the normalized error processes converge, or are at least tight, we say that the sequence $(u_n)$ is a rate, which, in addition, is sharp when the limiting process (or processes) is not trivial.

We suppose that $Y$ has no Gaussian part (otherwise a rate is known to be $u_n = \sqrt{n}$). Then rates are given in terms of the concentration of the Lévy measure of $Y$ around 0 and, further, we prove the convergence of the sequence $u_n(X^n - X)$ to a nontrivial limit under some further assumptions, which cover all stable processes and a lot of other Lévy processes whose Lévy measure behave like a stable Lévy measure near the origin. For example, when $Y$ is a symmetric stable process with index $\alpha \in (0,2)$, a sharp rate is $u_n = (n/\log n)^{1/\alpha}$; when $Y$ is stable but not symmetric, the rate is again $u_n = (n/\log n)^{1/\alpha}$ when $\alpha > 1$, but it becomes $u_n = n/(\log n)^2$ if $\alpha = 1$ and $u_n = n$ if $\alpha < 1$.

**1. Introduction.** We consider the following stochastic differential equation (SDE):

$$X_t = x_0 + \int_0^t f(X_{s-})\,dY_s, \tag{1.1}$$

where $f$ denotes a $C^3$ (three times differentiable) function and $Y$ is a Lévy process with characteristics $(b,c,F)$ with respect to the truncation function $h(x) = x\mathbb{1}_{\{|x|\le 1\}}$, that is,

$$E(e^{iuY_t}) = \exp t\left(iub - \frac{cu^2}{2} + \int F(dx)(e^{iux} - 1 - iux\mathbb{1}_{\{|x|\le 1\}})\right). \tag{1.2}$$

We also suppose that $f$ is such that (1.1) admits a (necessarily unique) *nonexploding* solution (this is the case, e.g., if $f$ has at most linear growth).

A number of papers have been devoted to studying the rate of convergence of the Euler scheme for this equation. That is, the approximated solution is defined at the times $i/n$, by induction on the integer $i$, according to the formula

$$X_0^n = x_0, \qquad X_{i/n}^n = X_{(i-1)/n}^n + f(X_{(i-1)/n}^n)(Y_{i/n} - Y_{(i-1)/n}). \tag{1.3}$$

This scheme allows for numerical computations, using Monte Carlo techniques, provided one can simulate the increments $Y_t - Y_s$ of the Lévy process $Y$: A first problem consists in computing an approximation of the expected value $E(h(X_1))$ for smooth enough functions $h$, and we need to evaluate the error $a_n(h) = E(h(X_1^n)) - E(h(X_1))$. A second problem is to compute an approximation of the law of some functional of the path, like, for example, $\sup_{t\leq 1} X_t$, and for this we need to evaluate the (discretized) error process, which is defined as

$$U_t^n = X_{[nt]/n}^n - X_{[nt]/n}. \tag{1.4}$$

Problem 1 has been extensively studied when $Y$ is continuous (i.e., $F = 0$) and $c > 0$: we can quote, with increasing order of generality as to the smoothness of $f$ and $h$, the works of Talay and Tubaro (1990) and Bally and Talay (1996a, b), where it is proved that $a_n(h)$ is of order $1/n$ and where an expansion of $a_n(h)$ as increasing powers of $1/n$ is even exhibited. In Protter and Talay (1997) the same problem is studied for discontinuous $Y$, but they only prove that $a_n(h) = O(1/n)$; see also a forthcoming paper by Kohatsu-Hida and Yoshida (2001) for an equation driven by a Wiener process plus a Poisson random measure. The techniques are essentially analytical.

For problem 2 one uses stochastic calculus techniques, and the idea is to find a *rate* $u_n$, that is, a sequence going to $\infty$ such that the sequence $(u_n U^n)$ is tight; the rate is called *sharp* if further the sequence $(u_n U^n)$ admits some limiting processes that are not identically 0. Even better is the case when the whole sequence $(u_n U^n)$ converges to a nondegenerate limit. In Jacod and Protter (1998) we have proved the following (more precise results are recalled below):

1. If $c > 0$, then a sharp rate is $u_n = \sqrt{n}$, and the sequence $(\sqrt{n} U^n)$ converges in law to a nondegenerate limit.
2. If $c = 0$ and $F$ is a finite measure, hence, $Y$ is a compound Poisson process plus a drift, then a sharp rate is $u_n = n$ if the drift $b$ is not 0; when $b = 0$, the rate is "infinite," meaning that for any $t$, we have $U_s^n = 0$ for all $s \leq t$ for $n$ large enough.
3. If $c = 0$ and $F$ is an infinite measure, then a rate is $u_n = \sqrt{n}$, but this rate is not sharp in the sense that $(\sqrt{n} U^n)$ goes in law to 0.

Although the implicit assumption that the increments of $Y$ can be simulated is somewhat unrealistic except in particular situations, which, however, include the case where $Y$ is a stable process, finding the exact rate of convergence is at least of much theoretical importance. Here we aim to find sharp rates for problem 2, when $c=0$ and $F(\mathbb{R})=\infty$. The crucial factor is the behavior of the Lévy measure $F$ near 0 (i.e., how many "small jumps" we have), which will be expressed through the following functions on $\mathbb{R}_+$:

$$\begin{aligned}\theta_+(\beta) &= F((\beta,\infty)),\\ \theta_-(\beta) &= F((-\infty,-\beta)), \qquad \theta(\beta)=\theta_+(\beta)+\theta_-(\beta).\end{aligned} \tag{1.5}$$

We introduce several assumptions, in which $\alpha$ denotes our basic index; here and below $C$ denotes a constant which may change from line to line, and may depend on $F$ just here, and also on $b$ and $f$ further below:

Hypothesis (H1-$\alpha$). *We have* $\theta(\beta)\le \frac{C}{\beta^\alpha}$ *for all* $\beta\in(0,1]$.

Hypothesis (H2-$\alpha$). *We have* $\beta^\alpha\theta_+(\beta)\to\theta_+$ *and* $\beta^\alpha\theta_-(\beta)\to\theta_-$ *as* $\beta\to 0$ *for some constants* $\theta_+,\theta_-\ge 0$, *and further*, $\theta:=\theta_++\theta_->0$. *We also set* $\theta'=\theta_+-\theta_-$, *and we observe that* $\theta(\beta)\sim\frac{\theta}{\beta^\alpha}$ *as* $\beta\to 0$.

Hypothesis (H3). *The measure $F$ is symmetrical about* 0.

Hypothesis (H4). *We have* $b=0$.

Note that Hypothesis (H2-$\alpha$) $\Rightarrow$ Hypothesis (H1-$\alpha$), and that Hypothesis (H1-2) always holds because $F$ integrates $x\mapsto|x|^2\wedge 1$, and Hypothesis (H1-0) [i.e., (H1-$\alpha$) for $\alpha=0$] holds iff the measure $F$ is finite, a case which we exclude. Under Hypothesis (H3) we have Hypothesis (H2-$\alpha$) as soon as $\theta(\beta)\sim\frac{\theta}{\beta^\alpha}$ as $\beta\to 0$, and $\theta_+=\theta_-=\theta/2$.

Unfortunately, we cannot totally fulfill our aim. But we find rates $u_n$ that are bigger than $\sqrt{n}$. And we prove that these rates are sharp and even that $u_nU^n$ converges in some reasonably general circumstances. Let us single out five different cases:

*Case* 1. We have Hypothesis (H1-$\alpha$) for some $\alpha>1$; then $u_n=(\frac{n}{\log n})^{1/\alpha}$.

*Case* 2a. We have Hypothesis (H1-$\alpha$) for $\alpha=1$; then $u_n=\frac{n}{(\log n)^2}$.

*Case* 2b. We have Hypothesis (H1-$\alpha$) for $\alpha=1$ and Hypothesis (H3); then $u_n=\frac{n}{\log n}$.

*Case* 3a. We have Hypothesis (H1-$\alpha$) for some $\alpha<1$; then $u_n=n$.

*Case* 3b. We have Hypothesis (H1-$\alpha$) for some $\alpha<1$ and Hypotheses (H3) and (H4); then $u_n=(\frac{n}{\log n})^{1/\alpha}$.

Clearly, Hypothesis (H1-$\alpha$) $\Rightarrow$ Hypothesis (H1-$\alpha'$) if $\alpha < \alpha'$, while the rate is better (i.e., bigger) when $\alpha$ decreases: one should take the smallest possible $\alpha$ for which Hypothesis (H1-$\alpha$) holds, although, of course, there might not be such a minimal $\alpha$. Observe also that the rate in Case 2b (resp. 3b) is strictly bigger than in Case 2a (resp. 3a): the symmetry of the driving process improves the quality of the Euler scheme under Hypothesis (H1-$\alpha$) when $\alpha \le 1$, while it does not affect the rate when $\alpha > 1$.

Now we describe the results of this paper. The first one concerns tightness [the assumption of $f$ is always that it is $C^3$ and that (1.1) has a nonexploding solution; this is not repeated in the next theorems].

THEOREM 1.1. *Assume that* $c = 0$ *and that Hypothesis* (H1-$\alpha$) *holds for some* $\alpha \in (0,2)$. *Then, with the above choice of* $u_n$, *the sequence* $(u_n U^n)$, *is tight.*

The results about limits necessitate the stronger Hypothesis (H2-$\alpha$) instead of Hypothesis (H1-$\alpha$), except in Case 3a; in all cases except 2a, the description of the limit invloves another process or additional random variables which are independent of $Y$, so we might need to enlarge the probability space to accomodate these.

Below, $\overline{Y}^n$ stands for the discretized process associated with $Y$, that is, $\overline{Y}^n_t = Y_{[nt]/n}$.

THEOREM 1.2. *Assume that* $c = 0$ *and that Hypothesis* (H1-$\alpha$) *holds for some* $\alpha \in (0,2)$. *Then in the following cases and with* $u_n$ *as above, the sequence* $(\overline{Y}^n, u_n U^n)$ *converges in law* (*for the Skorokhod topology*) *to* $(Y, U)$, *where* $U$ *is the unique solution* $U$ *of the linear equation*

$$U_t = \int_0^t f'(X_{s-}) U_{s-}\, dY_s - W_t, \tag{1.6}$$

*and where the process* $W$ *may be described as follows*:

(a) *In Case* 1, *and if further Hypothesis* (H2-$\alpha$) *holds*, *then*

$$W_t = \int_0^t f(X_{s-}) f'(X_{s-})\, dV_s, \tag{1.7}$$

*where* $V$ *is another Lévy process, independent of* $Y$ *and characterized by*

$$\begin{aligned} E(e^{iuV_t}) = \exp \frac{t\alpha}{2} \int & ((\theta_+^2 + \theta_-^2) \mathbb{1}_{\{x>0\}} + 2\theta_+\theta_- \mathbb{1}_{\{x<0\}}) \\ & \times \frac{1}{|x|^{1+\alpha}} (e^{iux} - 1 - iux)\, dx \end{aligned} \tag{1.8}$$

(*hence*, $V$ *is a stable process with index* $\alpha$).

(b) *In Case* 2a, *and if further Hypothesis* (H2-$\alpha$) *holds for* $\alpha = 1$, *then*

$$W_t = -\frac{(\theta_+ - \theta_-)^2}{4} \int_0^t f(X_{s-}) f'(X_{s-})\, ds, \tag{1.9}$$

*and we even have that* $u_n U^n$ *converges to* $U$ *in probability* (*locally uniformly in time*).

(c) *In Cases* 2b *and* 3b, *and if further Hypothesis* (H2-$\alpha$) *holds, then we have* (1.7), *where* $V$ *is another Lévy process, independent of* $Y$ *and characterized by*

$$E(e^{iuV_t}) = \exp t \int \frac{\theta^2 \alpha}{4} \frac{1}{|x|^{1+\alpha}} (e^{iux} - 1 - iux \mathbb{1}_{\{|x| \le 1\}})\, dx \tag{1.10}$$

(*hence,* $V$ *is a symmetric stable process with index* $\alpha$).

(d) *In Case* 3a, *then*

$$\begin{aligned} W_t = d \sum_{n:\, R_n \le t} ([f(X_{R_n}) - f(X_{R_n-})]\xi_n + f'(X_{R_n-})\Delta X_{R_n}(1-\xi_n)) \\ + \frac{d^2}{2} \int_0^t f(X_{s-}) f'(X_{s-})\, ds, \end{aligned} \tag{1.11}$$

*where* $d = b - \int_{\{|x| \le 1\}} xF(dx)$ *and* $(\xi_n)_{n \ge 1}$ *is a sequence of i.i.d. variables, uniform on* $[0,1]$ *and independent of* $Y$, *and* $(R_n)_{n \ge 1}$ *is an enumeration of the jump times of* $Y$ (*or of* $X$).

REMARK 1.1. For comparison with the cases excluded here and studied in Jacod and Protter (1998), let us mention that if $c = 0$ and $F$ is a finite measure [i.e., Hypothesis (H1-0) holds], then Theorem 1.2(d) holds without change. When $c > 0$, the sequence $(\overline{Y}^n, \sqrt{n} U^n)$ converges in law to $(Y, U)$, where $U$ solves (1.6) with

$$\begin{aligned} W_t = \sqrt{c} \sum_{n:\, R_n \le t} ([f(X_{R_n}) - f(X_{R_n-})]\sqrt{\xi_n}\kappa_n + f'(X_{R_n-})\Delta X_{R_n}\sqrt{1-\xi_n}\kappa'_n) \\ + \frac{c}{\sqrt{2}} \int_0^t f(X_{s-}) f'(X_{s-})\, dB_s \end{aligned}$$

and where $B$ is a standard Brownian motion, and $\xi_n$ is uniform over $[0,1]$, and $\kappa_n$ and $\kappa'_n$ are standard normal variables, all these being independent one from the other and from $Y$ as well.

REMARK 1.2. When $\theta_+ = \theta_-$ [e.g., under Hypothesis (H3)] then (1.8) and (1.10) agree (but, of course, for different values of $\alpha$). In Theorem 1.2(b) [resp. (d)], if $\theta' = \theta_+ - \theta_- = 0$ (resp. $d = 0$), the limiting process $U$ is identically 0. So these results are interesting only when $\theta' \ne 0$ (resp. $d \ne 0$), implying that $Y$ is dissymmetric, and otherwise the rate is not sharp.

REMARK 1.3. It would be possible, at the price of even more complicated computations, to accomodate other forms for Hypothesis (H2-$\alpha$): for example, if $\theta_+(\beta)$ and $\theta_-(\beta)$ are of order $\beta^{-\alpha}(\log\frac{1}{\beta})^\gamma$ as $\beta\to 0$ for some $\alpha\in(0,2)$ and $\gamma\in\mathbb{R}$. On the contrary, it seems rather difficult to express the rates $u_n$ directly in terms of the two functions $\theta_+(\beta)$ and $\theta_-(\beta)$.

REMARK 1.4. In Theorem 1.2(b) we have convergence in probability, so there ought to be an associated "central limit theorem." This suggests that we can improve the Euler scheme and simultaneously improve the rate, which is, indeed, the case: assume that Hypothesis (H1-$\alpha$) holds for $\alpha=1$. Then replace the Euler scheme of (1.3) by a modified Euler scheme $X'^n$ constructed as follows:

$$\begin{aligned} X'^n_0 &= x_0, \\ X'^n_{i/n} &= X'^n_{(i-1)/n} + f(X'^n_{(i-1)/n})(Y_{i/n}-Y_{(i-1)/n}) \\ &\qquad - f(X'^n_{(i-1)/n})f'(X'^n_{(i-1)/n})\gamma_n, \end{aligned} \tag{1.12}$$

where

$$\gamma_n = \frac{1}{2n^2}\int_{\{\log n/n<|x|\le 1\}} xF(dx)\int_{\{\log n/(n|x|)<y\le 1\}} yF(dy). \tag{1.13}$$

Denote by $U'^n$ the associated error process, that is, $U'^n_t = X'^n_{[nt]/n} - X_{[nt]/n}$. Then one can show that the sequence $(u'_n U'^n)$ is tight with $u'_n = n/\log n$. If further Hypothesis (H2-$\alpha$) holds with $\alpha=1$, then it is quite likely (although we have not proved it) that this sequence even converges in law to a nontrivial limiting process. The improvement in the rate, going from $n/(\log n)^2$ to $n/\log n$, is, of course, negligible in practice [not to speak about the fact that actually computing $\gamma_n$ in (1.13) might be a difficult task], so these results will not be proved below. Observe also that when further Hypothesis (H3) holds, then $\gamma_n$ above vanishes, so $X'^n = X^n$ and we recover Theorem 1.1 in Case 2b.

REMARK 1.5. We will introduce below—and use in a crucial way—a condition called (UT) [or sometimes (P-UT), for "predictably uniformly tight"] on a sequence of processes. Then one easily deduces from the proof of Theorem 1.1 that the sequence $(u_n U^n)$ satisfies the (UT) property, in addition to being tight.

REMARK 1.6. Let us say a word about our assumption that $f$ is $C^3$: it is fully used here for Case 1 only. For Cases 2 and 3 [statements (b)–(d) of Theorem 1.2 it is enough that $f$ be $C^2$, and an application of Theorem 3.5 of Kurtz and Protter (1991a, b), plus some scaling property, would give that the results hold as soon as $f$ is $C^1$, when $Y$ is a symmetric stable process.

Finally, let us mention that, for the sake of notational simplicity, we have considered only the one-dimensional case, but everything goes through in the multi-dimensional case as well, with exactly the same proofs.

The paper is organized as follows: in Section 2 we present a number of general tools connected with Euler approximations and limit theorem. More specific tools are developped in Section 3, while the heart of the proof (a long string of inequalities, extremely technical) is given in Section 4. The proofs of the two preceding theorems are given in Section 5.

**2. Some general tools.** Throughout all the paper, we suppose that $c=0$ and at least that Hypothesis (H1-$\alpha$) holds for some $\alpha\in(0,2)$. Recall also that $f$ is $C^3$ and (1.1) has a nonexploding solution.

2.1. *A condition for convergence of Euler schemes.* Let us write $Y'=f(X)$. We have

$$U_t^n=\sum_{i=1}^{[nt]}\int_{(i-1)/n}^{i/n}(f(X_{(i-1)/n}^n)-f(X_{(i-1)/n}))\,dY_s$$
$$-\sum_{i=1}^{[nt]}\int_{(i-1)/n}^{i/n}(Y'_{s-}-Y'_{(i-1)/n})\,dY_s.$$

Recalling $\overline{Y}_t^n=Y_{[nt]/n}$ and setting accordingly $\overline{X}^n=X_{[nt]/n}$ and

$$W_t^n=\sum_{i=1}^{[nt]}\int_{(i-1)/n}^{i/n}(Y'_{s-}-Y'_{(i-1)/n})\,dY_s, \tag{2.1}$$

we obtain

$$U_t^n=\int_0^t(f(\overline{X}_{s-}^n+U_{s-}^n)-f(\overline{X}_{s-}^n))\,d\overline{Y}_s^n-W_t^n. \tag{2.2}$$

Observe that the sequence $(\overline{Y}^n,\overline{X}^n)$ converges pointwise to $(Y,X)$ for the Skorohod topology. We can say more, and for this we recall the property (UT) defined in Jakubowski, Mémin and Pagès (1989) [see also Jacod and Shiryaev (2003), Chapter VI.6, where it is called (P-UT)]. Let $Z^n$ be a sequence of semimartingale, with the canonical decompositions

$$Z_t^n=A_t^{n,a}+M_t^{n,a}+\sum_{s\le t}\Delta Z_s^n 1_{\{|\Delta Z_s^n|>a\}}, \tag{2.3}$$

where $a>0$ and $A^{n,a}$ is a predictable process with locally bounded variation and $M^{n,a}$ is a (locally bounded) local martingale. Then we say that the sequence $(Z^n)$ satisfies (UT) if for any $t<\infty$, the sequence of real-valued random variables

$$\mathrm{Var}(A^{n,a})_t + \langle M^{n,a}, M^{n,a}\rangle_t + \sum_{s\le t} |\Delta Z_s^n| \mathbb{1}_{\{|\Delta Z_s^n|>a\}}$$

is tight. This property does not depend on the choice of $a\in(0,\infty)$.

The following lemma applies, in particular, when $\Gamma^n = \overline{Y}^n$ because then $\Gamma_1^n = Y_1$. Its setting is as follows: we have a triangular array of rowwise i.i.d. $d$-dimensional random variables $\zeta_i^n, i=1,2,\dots,$ and we set $\Gamma_t^n = \sum_{i=1}^{[nt]} \zeta_i^n$.

LEMMA 2.1. *If $\Gamma_1^n$ converges in law to a limit $U$, then there is a $d$-dimensional Lévy process $\Gamma$ such that $\Gamma_1 = U$; this process is unique in law and $\Gamma^n$ converges in law to $\Gamma$ (for the Skorokhod topology). Further, the sequence $(\Gamma^n)$ has (UT).*

PROOF. The first claims are all well known [see, e.g., Jacod and Shiryaev (2003), Chapter VII.3.6], and only the last one needs proving. Since the (UT) property holds for a multi-dimensional sequence iff it holds separately for the sequence of each component, for the last claim we may assume w.l.o.g. that the variables are one-dimensional. Let $(b,c,F)$ be the characteristics of the Lévy process $\Gamma$, and take $a>1$ such that $F(\{x:|x|=a\}) = F(\{x:|x| = -a\}) = 0$. Set $b_n = E(\zeta_i^n \mathbb{1}_{\{|\zeta_i^n|\le a\}})$ and $c_n = E((\zeta_i^n)^2 \mathbb{1}_{\{|\zeta_i^n|\le a\}}) - b_n^2$ and $\gamma_n = P(|\zeta_i^n|>a)$, for all $a\in D$ we have, by virtue of Jacod and Shiryaev (2003),

$$(2.4)\qquad \begin{aligned} &nb_n \to b + \int_{\{1<|x|\le a\}} xF(dx),\\ &nc_n \to \int_{\{|x|\le a\}} x^2 F(dx), \qquad n\gamma_n \to F(x:|x|>a). \end{aligned}$$

Writing (2.3) for $\Gamma^n$ gives $A_t^{n,a} = b_n[nt]$ and $\langle M^{n,a}, M^{n,a}\rangle_t = c_n[nt]$. Set further $V_t^{n,a} = \sum_{s\le t} |\Delta\Gamma_s^n| \mathbb{1}_{\{|\Delta\Gamma_s^n|>a\}}$ and $V_t'^{n,a} = \sum_{s\le t} \mathbb{1}_{\{|\Delta\Gamma_s^n|>a\}}$ and $H_t^n = \sup_{s\le t} |\Gamma_s^n|$. Then (2.4) yields that the sequences $(A_t^{n,a})_n$, $(\langle M^{n,a}, M^{n,a}\rangle_t)_n$ and $(V_t'^{n,a})_n$ are tight (for the later, note that $E(V_t'^{n,a}) \le [nt]\gamma_n$), and the convergence in law $\Gamma^n\to\Gamma$ yields the tightness of the sequence $(H_t^n)_n$ for all $t$. Since $V_t^{n,a} \le H_t^n V_t'^{n,a}$, the result is obvious. □

Using (2.2) and Lemma 2.1 and the fact that $(\overline{Y}^n, \overline{X}^n)$ converges to $(Y,X)$, as mentioned earlier, and following the proof of Theorem 3.2 of Jacod and Protter (1998), which itself is based upon Kurtz and Protter (1991b, 1996) [see also Słomiński (1989) and Mémin and Słomiński (1991)], we readily obtain the following:

THEOREM 2.2. *Let $(u_n)$ be a sequence of reals increasing to $+\infty$, such that the sequence $(\overline{Y}^n, u_n W^n)$ converges to a limit $(Y, W)$ in law, where $W$ is possibly defined on an extension of the original probability space (resp. in probability, with $W$ given on the original space). Then the sequence $(\overline{Y}^n, u_n U^n)$ converges in law (resp. in probability) to $(Y, U)$, where $U$ is the unique solution $U$ of the following equation:*

$$U_t = \int_0^t f'(X_{s-}) U_{s-}\, dY_s - W_t. \tag{2.5}$$

Up to taking subsequences, we deduce the following:

COROLLARY 2.3. *Let $(u_n)$ be a sequence of reals increasing to $+\infty$, such that the sequence $(\overline{Y}^n, u_n W^n)$ is tight. Then the sequence $(u_n U^n)$ is also tight.*

2.2. *Localization.* In order to avoid a lot of technical problems, we will "localize" it in the sense of the following proposition:

PROPOSITION 2.4. *Suppose that either Theorem 1.1 or 1.2 holds for all Lévy processes $Y$ satisfying the relevant assumptions and having further bounded jumps, and for all $C^3$ functions $f$ with compact support. Then these theorems also hold for a Lévy process $Y$ with unbounded jumps and a $C^3$ function $f$ with noncompact support.*

PROOF. We start with a Lévy process $Y$ satisfying the assumptions of one of our theorems and with a $C^3$ function $f$. We call $W$ the limit of $u_n W^n$ in case of Theorem 1.2.

For $p \in \mathbb{N}^\star$ we consider the new Lévy process $Y(p)_t = Y_t - \sum_{s \le t} \Delta Y_s \mathbb{1}_{\{|\Delta Y_s| > p\}}$, and a $C^3$ function $f_p$ with compact support, satisfying $f_p(x) = f(x)$ for $|x| \le p$. Finally, we associate with $Y(p)$ and $f_p$ the same terms with $Y$ and $f$, writing, in particular $X(p)$, $\overline{Y}(p)^n$, $U(p)^n$ instead of $X$, $\overline{Y}^n$, $U^n$.

Observe that $Y(p)$ satisfies the same Hypotheses (H1-$\alpha$) or (H2-$\alpha$) or (H3) or (H4) than $Y$, and the numbers $\theta_+$ and $\theta_-$ in Hypothesis (H2-$\alpha$) or $d$ in Theorem 1.2(d) are the same for $Y$ and each $Y(p)$. So our hypothesis yields in case of Theorem 1.1 that for each $p$ the sequence $(u_n U(p)^n)$ is tight. For Theorem 1.2 it yields that for each $p$ the sequence $(\overline{Y}(p)^n, u_n U(p)^n)$ converges to $(Y(p), U(p))$, where $U(p)$ satisfies (1.6) relative to some process $W(p)$: for (a) or (c) this process is given by (1.7) with $V(p) = V$ independent of $p$ and with $X(p)$; for (b) it is given by (1.9) with $X(p)$; for (d) we have to be more careful: for each $q \in \mathbb{N}^\star$, we denote by $(R(q)_n)_{n \ge 1}$ an enumeration of all jump times of $Y$ with size in $(q-1, q]$, and let $(\xi(q)_n)_{n,q \ge 1}$ be a double sequence of i.i.d. variables uniform over $[0,1]$ and independent of $Y$; then $W(p)$ can be taken to be the process defined by (1.11), where the summation in the first term extends to all $R(q)_n$ and $\xi(q)_n$ for $q \le p$ and $X$ is replaced by $X(p)$; similarly, for $W$ we can take $W = W(\infty)$, with the process $X$.

Set $S_p = \inf(t : |X_t| \ge p$ or $|\Delta Y_t| > p)$. We have $Y = Y(p)$, hence, $X = X(p)$ and $U^n = U(p)^n$ and $\overline{Y}^n = \overline{Y}(p)^n$, on the interval $[0, S_p)$; in case of Theorem 1.2, we also have $W = W(p)$, hence, $U = U(p)$ as well, on the interval $[0, S_p)$.

Let us first consider the case of Theorem 1.1. For any $\varepsilon > 0$ and $t > 0$, there is a $p$ such that $P(S_p \le t) \le \varepsilon$, and a compact set $K$ in the Skorokhod space $\mathbb{D}([0,\infty), \mathbb{R})$ which depends on the sample path only up to time $t$ and such that $P(u_n U(p)^n \notin K) \le \varepsilon$. Since $U = U(p)$ on $[0, S_p)$, we have $\{u_n U^n \notin K, S_p > t\} = \{u_n U(p)^n \notin K, S_p > t\}$, hence, $P(u_n U^n \notin K) \le 2\varepsilon$, and this proves the tightness of the sequence $(u_n U^n)$.

Let us next consider the case of Theorem 1.2. For any continuous bounded function $\Phi_t$ on the Skorokhod space $\mathbb{D}([0,\infty), \mathbb{R}^2)$ which depends on the sample path only up to time $t$, we have

$$|E(\Phi_t(\overline{Y}^n, u_n U^n)) - E(\Phi_t(\overline{Y}(p)^n, u_n U(p)^n))| \le 2\|\Phi_t\| P(S_p \le t),$$

$$|E(\Phi_t(Y, U)) - E(\Phi_t(Y(p), U(p)))| \le 2\|\Phi_t\| P(S_p \le t).$$

Since $P(S_p \le t) \to 0$ and since $E(\Phi_t((\overline{Y}(p)^n, u_n U(p)^n)) \to E(\Phi_t((Y(p), U(p)))$ for every $t$ as $p \to \infty$, we get $E(\Phi_t((\overline{Y}^n, u_n U^n)) \to E(\Phi_t((Y(p), U(p)))$, hence, the result for the convergence in law. For the convergence in measure in Theorem 1.2(b), the proof is similar. □

2.3. *Some limit theorems.* In a rather natural way, the solution to our problem goes through various limit theorem concerning sums of triangular arrays of the form

$$\Gamma_t^n = \sum_{i=1}^{[nt]} \zeta_i^n, \tag{2.6}$$

where for each $n$ we have $\mathbb{R}^d$-valued random variables $(\zeta_n^n)_{i \ge 1}$ such that each $\zeta_i^n$ is $\mathcal{F}_{i/n}$-measurable. Below we give various conditions (very far from being optimal) insuring tightness or convergence of the sequence $(\Gamma^n)$.

First we introduce a set of conditions, where $\xi_n$, $\xi'_n$, $\xi''_{n,y}$ denote arbitrary finite constants:

$$E(|\zeta_i^n| | \mathcal{F}_{(i-1)/n}) \le \frac{\xi_n}{n}, \tag{2.7}$$

$$\begin{aligned} |E(\zeta_i^n | \mathcal{F}_{(i-1)/n})| &\le \frac{\xi_n}{n}, \\ E(|\zeta_i^n|^2 | \mathcal{F}_{(i-1)/n}) &\le \frac{\xi'_n}{n}, \end{aligned} \tag{2.8}$$

$$
\begin{aligned}
&|E(\zeta_i^n 1_{\{|\zeta_i^n|\le 1\}}|\mathcal{F}_{(i-1)/n})| \le \frac{\xi_n}{n},\\
&E(|\zeta_i^n|^2 1_{\{|\zeta_i^n|\le 1\}}|\mathcal{F}_{(i-1)/n}) \le \frac{\xi'_n}{n}, \qquad (2.9)\\
&P(|\zeta_i^n| > y|\mathcal{F}_{(i-1)/n}) \le \frac{\xi''_{n,y}}{n} \qquad \forall y > 1.
\end{aligned}
$$

Note that (2.8) with $\widehat{\xi}_n$ and $\widehat{\xi}'_n$ implies (2.9) with $\xi_n = \widehat{\xi}_n + \widehat{\xi}'_n$ and $\xi'_n = \widehat{\xi}'_n$ and $\xi''_{n,y} = \widehat{\xi}'_n/y^2$. Also, (2.7) with $\widehat{\xi}_n$ implies (2.9) with $\xi_n = \xi'_n = \widehat{\xi}_n$ and $\xi_{n,y} = \widehat{\xi}_n/y$. Observe also that if (2.9) holds, then the last inequality is true for $y \in (0,1]$ as well, with $\xi''_{n,y} = \xi''_{n,1} + \xi'_n/y^2$.

Part (a) below is well known, and part (b) follows from Theorem VI.5.10 of Jacod and Shiryaev (2003). By $\Gamma^n \xrightarrow{P} 0$, we mean that $\sup_{s\le t}|\Gamma_s^n|$ goes to 0 in probability for all $t$.

LEMMA 2.5. (a) *For $\Gamma^n \xrightarrow{P} 0$, it is enough that either (2.7) or (2.8) or (2.9) hold with*

$$
\lim_n \xi_n = 0, \qquad \lim_n \xi'_n = 0, \qquad \lim_n \xi''_{n,y} = 0 \qquad \forall y > 1. \qquad (2.10)
$$

(b) *For the sequence $(\Gamma^n)$ to be tight for the Skorokhod topology, it is enough that the sequence of each of the $d$ components of $\zeta_i^n$ satisfies either (2.7) or (2.8) or (2.9) with*

$$
\begin{aligned}
&\limsup_n \xi_n < \infty, \qquad \limsup_n \xi'_n < \infty,\\
&\limsup_n \xi''_{n,y} < \infty, \qquad \lim_{y\uparrow\infty}\limsup_n \xi''_{n,y} = 0.
\end{aligned} \qquad (2.11)
$$

The conditions in Lemma 2.5 can be substituted with conditions on the following conditional characteristic functions:

LEMMA 2.6. *Suppose that one can find constants $\xi'''_{n,v}$ such that*

$$
\sup_{u:|u|\le v} |1 - E(e^{iu.\zeta_i^n}|\mathcal{F}_{(i-1)/n})| \le \frac{\xi'''_{n,v}}{n}, \qquad (2.12)
$$

*then (2.9) holds with $\xi_n = \xi'_n = C\xi'''_{n,1}$ and $\xi''_{n,y} = C\xi'''_{n,1/y}$.*

PROOF. It is enough to consider the one-dimensional case $d = 1$. We use known facts about characteristic functions, which readily pass to "conditional" characteristic functions. We have [see, e.g., (2) in the proof of Lemma VII.2.16 of Jacod and Shiryaev (2003)]

$$
E(|w\zeta_i^n|^2 \wedge 1|\mathcal{F}_{(i-1)/n}) \le C\int_{\{|u|\le w\}} |1 - E(e^{iu.\zeta_i^n}|\mathcal{F}_{(i-1)/n})|\,du \le \frac{C}{n}\xi'''_{n,w}.
$$

This readily gives (2.9) with $\xi'_n = C\xi'''_{n,1}$ and $\xi''_{n,y} = C\xi'''_{n,1/y}$. We also have the estimate $|x1_{\{|x|\le 1\}} - \sin x| \le x^2 \wedge 1$, hence,

$$|E(\zeta_i^n 1_{\{|\zeta_i^n|\le 1\}}|\mathcal{F}_{(i-1)/n})| \le |1 - E(e^{i\zeta_i^n}|\mathcal{F}_{(i-1)/n})| + E(|\zeta_i^n|^2 \wedge 1|\mathcal{F}_{(i-1)/n}),$$

and (2.9) holds with $\xi_n = C\xi'''_{n,1}$. □

LEMMA 2.7. *In the previous setting, suppose that we have* (2.9) *with* $\widehat{\xi}_n$, $\widehat{\xi}'_n$ *and* $\widehat{\xi}''_{n,y}$ *and* $\widehat{C}_y$, *and assume that* $\widehat{\xi}_n/n \to 0$. *Then the variables* $\zeta_i'^n = \zeta_i^n - E(\zeta_i^n 1_{\{\zeta_i^n|\le 1\}}|\mathcal{F}_{(i-1)/n})$ *satisfy* (2.9) *with, for all* $n$ *large enough,*

$$\xi_n = 6\,\widehat{\xi}''_{n,1/2}, \qquad \xi'_n = 4\,\widehat{\xi}'_n + 8\,\widehat{\xi}''_{n,1} + \frac{4\widehat{\xi}_n^2}{n}, \qquad \xi''_{n,y} = \widehat{\xi}''_{n,y-1/2}. \tag{2.13}$$

PROOF. Set $a_i^n = E(\zeta_i^n 1_{\{|\zeta_i^n|\le 1\}}\mathcal{F}_{(i-1)/n})$. We have $|a_i^n| \le \widehat{\xi}_n/n$, so $\varepsilon_n = \sup_i |a_i^n| \to 0$ and, up to taking $n$ large enough, we can assume that $\varepsilon_n \le 1/2$. Then

$$\begin{aligned} y > 1 \quad \Longrightarrow \quad \{|\zeta_i'^n| > y\} &\subset \{|\zeta_i^n| > (y - 1/2)\}, \\ (\zeta_i'^n)^2 1_{\{|\zeta_i'^n|\le 1\}} &\le 4((\zeta_i^n)^2 1_{\{|\zeta_i^n|\le 1\}} + \varepsilon_n^2 + 2 1_{\{|\zeta_i'^n|>1\}}), \\ |\zeta_i'^n 1_{\{|\zeta_i'^n|\le 1\}} - (\zeta_i^n 1_{\{|\zeta_i^n|\le 1\}} - a_i^n)| &\le 3(1_{\{|\zeta_i'^n|>1\}} + 1_{\{|\zeta_i^n|>1\}}). \end{aligned}$$

The result is then obvious. □

Finally, we will often encounter the following situation, in connection with our basic process $Y$: we have a pair $(Z^n, \Gamma^n)$ of (possibly multi-dimensional) processes of the form

$$Z_t^n = \sum_{i=1}^{[nt]} \eta_i^n, \qquad \Gamma_t^n = \sum_{i=1}^{[nt]} \zeta_i^n. \tag{2.14}$$

Further, we have

$$\zeta_i^n = g(X_{(i-1)/n})\zeta_i'^n, \tag{2.15}$$

and for each $n$ the sequence $(\eta_i^n, \zeta_i'^n), i = 1, 2, \ldots,$ is i.i.d. We set $\Gamma_t'^n = \sum_{i=1}^{[nt]} \zeta_i'^n$.

Assume also that $Z^n$ converges in probability (for the Skorokhod topology) to a limit $Z$ of the form $Z_t = Y_t + at$ for some constant $a$. Then, combining Lemma 2.1 with a fundamental property of convergent sequences of processes having (UT), we get the following:

LEMMA 2.8. *In the previous setting, suppose that the pair* $(Z_1^n, \Gamma_1'^n)$ *of random variables converges in law to* $(Z_1, \gamma')$ *with* $\gamma'$ *a random variable independent of* $Z_1$*, and that* $g$ *is a continuous function. Then there is a Lévy process* $\Gamma'$*, independent of* $Y$ *and unique in law, such that the processes* $(Z^n, \Gamma^n, \Gamma'^n)$ *converge in law to* $(Z, \Gamma, \Gamma')$*, where* $\Gamma_t = \int_0^t g(X_{s-})\, d\Gamma'_s$*. If further* $\gamma'$ *is a constant, then we get* $\Gamma_t = \int_0^t g(X_{s-})\gamma'\, ds$*, and the convergence of* $(Z^n, \Gamma^n, \Gamma'^n)$ *takes place in probability.*

PROOF. Lemma 2.1 yields the convergence in law of $(Z^n, \Gamma'^n)$ to some Lévy process $(Z', \Gamma')$ which is unique in law, and the independence of the variables $Z_1$ and $\gamma'$ implies the independence of the processes $Z$ and $\Gamma'$. Further, since $Z'_1 = Z_1$ in distribution, then the laws of $Z$ and $Z'$ are the same: so we can realize $(Z', \Gamma')$ with a first component $Z'$ equal to the original process $Z$. And if further $\gamma'$ is a constant, then obviously $\Gamma'_t = \gamma' t$ and the convergence of $(Z^n, \Gamma'^n)$ to $(Z, \Gamma')$ holds in probability.

Finally, Lemma 2.1 also yields the (UT) property for $(\Gamma'^n)$, and a fundamental property of (UT) [Theorem VI.6.22 of Jacod and Shiryaev (2003)] gives the claim. □

**3. Some preliminaries.** As seen in Proposition 2.4, we can and will assume that $f$ is $C^3$ with compact support, and that $|\Delta Y| \le p$ identically for some integer $p \ge 1$, which amounts to saying that $\theta(p) = 0$.

3.1. *About the Lévy measure.* The following quantities, where $\beta > 0$, will be of interest:

$$
\begin{aligned}
&c(\beta) = \int_{\{|x|\le\beta\}} |x|^2 F(dx), \\
&d_+(\beta) = \int_{\{x>\beta\}} |x| F(dx), \qquad d_-(\beta) = \int_{\{x<-\beta\}} |x| F(dx), \\
&\rho_+(\beta) = \int_{\{x>\beta\}} |x|^\alpha F(dx), \qquad \rho_-(\beta) = \int_{\{x<-\beta\}} |x|^\alpha F(dx), \\
&\delta(\beta) = d_+(\beta) + d_-(\beta), \qquad \rho(\beta) = \rho_+(\beta) + \rho_-(\beta), \\
&d'(\beta) = d_+(\beta) - d_-(\beta), \qquad b' = b + \int_{\{|x|>1\}} x F(dx), \\
&d(\beta) = b' - d'(\beta).
\end{aligned}
\tag{3.1}
$$

Note that $d(\beta) = b - \int_{\{\beta<|x|\le 1\}} x F(dx)$ if $\beta < 1$. We will now give some estimates on these quantities. First, observe that, for all $0 \le a < b \le 1$ and

$\gamma > 0$, we have

$$\int_{\{a<x\le b\}} |x|^\gamma F(dx) = \gamma \int_0^b y^{\gamma-1}(\theta_+(y\vee a) - \theta_+(b))\,dy, \tag{3.2}$$

and a similar relation on the negative side. Introduce also the notation

$$s(\beta) = \begin{cases} 1, & \text{if } \alpha < 1, \\ \log\left(\dfrac{1}{\beta}\right), & \text{if } \alpha = 1, \\ \dfrac{1}{\beta^{\alpha-1}}, & \text{if } \alpha > 1. \end{cases}$$

Then, since $\theta(p) = 0$, we readily deduce that under Hypothesis (H1-$\alpha$) we have

$$\begin{gathered} c(\beta) \le C\beta^{2-\alpha}, \qquad \rho(\beta) \le C\log\left(\frac{1}{\beta}\right), \qquad \int_{\{|x|>\beta\}} |x|^{\alpha/2} F(dx) \le \frac{C}{\beta^{\alpha/2}}, \\ \delta(\beta) + |d(\beta)| + d_+(\beta) + d_-(\beta) + |d'(\beta)| \le Cs(\beta). \end{gathered} \tag{3.3}$$

Further, if $\alpha < 1$, then $d_+(\beta)$, $d_-(\beta)$ and $d(\beta)$ converge as $\beta \to 0$ to some finite limits $d_+$, $d_-$ and $d$, and $d$ is as in Theorem 1.2(d).

Under Hypothesis (H3) we also have

$$d'(\beta) = 0, \qquad |d(\beta)| \le C, \tag{3.4}$$

while under Hypotheses (H3) and (H4) we even have $b = b' = d(\beta) = 0$.

Next suppose that Hypothesis (H2-$\alpha$) holds. Taking advantage of (3.2), we obtain the following equivalences or convergences as $\beta \to 0$:

$$\begin{aligned} & c(\beta) \sim \frac{\alpha\theta}{2-\alpha}\beta^{2-\alpha}, \\ & \frac{\rho_+(\beta)}{\log(1/\beta)} \to \alpha\theta_+, \qquad \frac{\rho_-(\beta)}{\log(1/\beta)} \to \alpha\theta_-, \\ & d_+(\beta) \to d_+, \qquad d_-(\beta) \to d_- \qquad \text{if } \alpha < 1, \\ & \frac{d_+(\beta)}{\log(1/\beta)} \to \theta_+, \qquad \frac{d_-(\beta)}{\log(1/\beta)} \to \theta_- \qquad \text{if } \alpha = 1, \\ & \beta^{\alpha-1} d_+(\beta) \to \frac{\alpha\theta_+}{\alpha-1}, \qquad \beta^{\alpha-1} d_-(\beta) \to \frac{\alpha\theta_-}{\alpha-1} \qquad \text{if } \alpha > 1. \end{aligned} \tag{3.5}$$

We will also need an estimate of the integral of $x\log x$ w.r.t. $F$ when $\alpha = 1$. For this we first observe that, analogously to (3.2),

$$\int_{\{a<x\le b\}} (x\log x) F(dx) = \int_0^b (1+\log y)(\theta_+(y\vee a) - \theta_+(b))\,dy,$$

and a similar relation on the negative side. We then deduce that for every $a > 0$ and as $\beta \to 0$,

$$\text{(3.6)}\quad \text{Hypothesis (H2-1)} \Rightarrow \begin{cases} \dfrac{1}{(\log(1/\beta))^2}\displaystyle\int_\beta^a (x\log x)F(dx) \to -\dfrac{\theta_+}{2}, \\ \dfrac{1}{(\log(1/\beta))^2}\displaystyle\int_{-a}^{-\beta} (|x|\log|x|)F(dx) \to -\dfrac{\theta_-}{2}. \end{cases}$$

3.2. *About the Lévy process.* Now we split the processes $Y$ and $Y' = f(X)$. We first recall that if $\mu$ is the jump measure of $Y$ and $\nu(ds,dx) = ds \otimes F(dx)$ is its predictable compensator, for each $\beta > 0$, we can write [recalling that $c = 0$; we denote by $U * (\mu - \nu)$ the stochastic integral process of the predictable function $U$ on $\overline{\Omega} \times \mathbb{R}_+ \times \mathbb{R}$ w.r.t. $\mu - \nu$]

$$\text{(3.7)}\quad \begin{aligned} &Y = A^\beta + M^\beta + N^\beta \qquad \text{where} \\ &A^\beta_t = d(\beta)t, \qquad M^\beta = x\mathbb{1}_{\{|x|\le\beta\}} * (\mu - \nu), \qquad N^\beta = x\mathbb{1}_{\{|x|>\beta\}} * \mu. \end{aligned}$$

Then $M^\beta$ is a square-integrable martingale with predictable bracket $\langle M^\beta, M^\beta\rangle_t = c(\beta)t$. Since $|\Delta Y| \le p$, we have $N^p = 0$ and $Y = A^p + M^p$ and $A^p_t = b't$. We also have $A^1_t = bt$.

Next, we set $g = ff'$, which is a $C^2$ function with compact support. We have the decomposition

$$\text{(3.8)}\qquad G(x,y) := f(x + yf(x)) - f(x) = yg(x) + y^2k(x,y),$$

with $k$ a $C^1$ function which vanishes outside $K \times \mathbb{R}$ for some compact subset $K$ of $\mathbb{R}$.

Now we turn to the decomposition of the semimartingale $Y' = f(X)$. We introduce the notation

$$b^\beta_t = g(X_{t-})\,d(\beta) + \int_{\{|x|\le\beta\}} F(dx)x^2k(X_{t-},x).$$

Then Itô's formula yields

$$Y' = Y'_0 + A'^\beta + M'^\beta + N'^\beta,$$

where

$$\begin{aligned} A'^\beta_t &= \int_0^t b^\beta_s\,ds, \\ M'^\beta &= G(X_-,x)\mathbb{1}_{\{|x|\le\beta\}} * (\mu - \nu), \\ N'^\beta &= G(X_-,x)\mathbb{1}_{\{|x|>\beta\}} * \mu. \end{aligned}$$

Observe that

$$\text{(3.9)}\qquad |b^\beta_t| \le C(|d(\beta)|) + c(\beta)$$

and

$$\langle M'^\beta, M'^\beta\rangle_t = \int_0^t c'_{\beta,s}\,ds \qquad \text{where } c'_{\beta,t} \le Cc(\beta). \tag{3.10}$$

We also set

$$Y^\beta = A^\beta + M^\beta, \qquad Y'^\beta = A'^\beta + M'^\beta. \tag{3.11}$$

3.3. *A decomposition for $W^n$.* 1. The rates $u_n$ have been described in Section 1 to the case we are in. We will also choose a sequence $\beta_n$ going to 0 in such a way that

$$\lambda_n = \frac{\theta(\beta_n)}{n} \to 0. \tag{3.12}$$

We write $c_n = c(\beta_n)$, $d_n = d(\beta_n)$, $d'_n = d'(\beta_n)$, $\rho_n = \rho(\beta_n)$ and $\delta_n = \delta(\beta_n)$. The precise choice of $\beta_n$ is as follows (we repeat also the definition of $u_n$ for easier reading):

*Case* 1. $u_n = (\frac{n}{\log n})^{1/\alpha}$ and $\beta_n = \frac{\log n}{n^{1/(2\alpha)}}$.

*Case* 2a. $u_n = \frac{n}{(\log n)^2}$ and $\beta_n = \frac{\log n}{n}$.

*Case* 2b. $u_n = \frac{n}{\log n}$ and $\beta_n = \frac{\log n}{n}$.

*Case* 3a. $u_n = n$ and $\beta_n = \frac{(\log n)^2}{n}$.

*Case* 3b. $u_n = (\frac{n}{\log n})^{1/\alpha}$ and $\beta_n = (\frac{\log n}{n})^{1/\alpha}$.

Taking advantage of the estimates of (3.3) we get the following:

$$\begin{aligned}
&\textit{Case } 1. && c_n \le C\frac{(\log n)^{2-\alpha}}{n^{(2-\alpha)/(2\alpha)}}, \qquad |d_n| + \delta_n \le C\frac{n^{\alpha-1/(2\alpha)}}{(\log n)^{\alpha-1}},\\
& && \lambda_n \le C\frac{1}{n^{1/2}(\log n)^\alpha}\\
&\textit{Case } 2a. && c_n \le C\frac{\log n}{n}, \qquad |d_n| + \delta_n \le C\log n, \qquad \lambda_n \le \frac{C}{\log n},\\
&\textit{Case } 2b. && c_n \le C\frac{\log n}{n}, \qquad |d_n| \le C, \qquad \delta_n \le C\log n, \qquad \lambda_n \le \frac{C}{\log n},\\
&\textit{Case } 3a. && c_n \le C\frac{(\log n)^{4-2\alpha}}{n^{2-\alpha}}, \qquad |d_n| + \delta_n \le C,\\
& && \lambda_n \le \frac{C}{n^{1-\alpha}(\log n)^{2\alpha}},\\
&\textit{Case } 3b. && c_n \le C\left(\frac{\log n}{n}\right)^{(2-\alpha)/\alpha}, \quad d_n = 0, \qquad \delta_n \le C, \qquad \lambda_n \le \frac{C}{\log n}.
\end{aligned} \tag{3.13}$$

If, further, Hypothesis (H2-$\alpha$) holds, we get (since $d_n = b' - d'_n$) the following:

$$
\begin{aligned}
&\textit{Case } 1. && c_n \sim \frac{\alpha\theta}{2-\alpha}\frac{(\log n)^{2-\alpha}}{n^{(2-\alpha)/(2\alpha)}}, \qquad \frac{d_n(\log n)^{\alpha-1}}{n^{(\alpha-1)/(2\alpha)}} \to -\frac{\alpha\theta'}{\alpha-1},\\
& && \lambda_n \sim \theta\frac{1}{n^{1/2}(\log n)^{\alpha}},\\
&\textit{Case } 2\text{a}. && c_n \sim \theta\frac{\log n}{n}, \qquad \frac{d_n}{\log n} \to -\theta', \qquad \lambda_n \sim \frac{\theta}{\log n},\\
&\textit{Case } 3\text{a}. && c_n \sim \frac{\alpha\theta}{2-\alpha}\frac{(\log n)^{4-2\alpha}}{n^{2-\alpha}}, \qquad d_n \to d,\\
& && \lambda_n \sim \frac{\theta}{n^{1-\alpha}(\log n)^{2\alpha}},\\
&\textit{Case } 2\text{b}, 3\text{b}. && \lambda_n \sim \frac{\theta}{\log n}.
\end{aligned} \tag{3.14}
$$

2. By virtue of Theorem 2.2 we have to prove the convergence of the processes $(\overline{Y}^n, u_n W^n)$. Both $\overline{Y}^n_t$ and $u_n W^n_t$ are the sums for $i$ between 1 and $[nt]$ of i.i.d. variables, say $y^n_i = Y_{i/n} - Y_{(i-1)/n}$ and $w^n_i$, each one depending only on the increments of $Y$ over the interval $I(n,i) = (\frac{i-1}{n}, \frac{i}{n}]$, and for $w^n_i$ on the "truncation" at level $\beta_n$.

Each of these variables, $y^n_i$ and $w^n_i$, will in turn be decomposed into small bits which are handled separately, and which we call $\zeta^n_i(j)$ for $j$ between 1 and 14! This is done in such a way that if

$$\Gamma^n(j)_t = \sum_{i=1}^{[nt]} \zeta^n_i(j),$$

then

$$\overline{Y}^n = \Gamma^n(13) + \Gamma^n(14), \qquad u_n W^n = \sum_{j=1}^{12} \Gamma^n(j). \tag{3.15}$$

In order to do this, we introduce a number of notations. First, denote by $X'^{n,i}$ the unique solution over $(\frac{i-1}{n}, \infty)$ of the equation $dX'^{n,i} = f(X'^{n,i}_-)\,dY^{\beta_n}$ starting at $X_{(i-1)/n}$ at time $\frac{i-1}{n}$. Then we set

$$
\begin{aligned}
&Y^{n,i}_t = Y^{\beta_n}_t - Y^{\beta_n}_{(i-1)/n}, \qquad Y'^{n,i}_t = Y'^{\beta_n}_t - Y'^{\beta_n}_{(i-1)/n}\\
&M^{n,i}_t = M^{\beta_n}_t - M^{\beta_n}_{(i-1)/n}, \qquad A^{n,i}_t = A^{\beta_n}_t - A^{\beta_n}_{(i-1)/n} = d_n\left(t - \frac{i-1}{n}\right)\\
&\hfill \text{for } t \ge \frac{i-1}{n},
\end{aligned}
$$

$$\widehat{Y}_i^n = \sup_{s \in I(n,i)} |Y_s^{n,i}|, \qquad \widehat{M}_i^n = \sup_{s \in I(n,i)} |M_s^{n,i}|, \qquad \widehat{A}_i^n = \sup_{s \in I(n,i)} |A_s^{n,i}|,$$

$$\widetilde{X}_i^n = \sup_{s \in I(n,i)} |X_s - X_{(i-1)/n}|, \qquad \widetilde{X}_i'^n = \sup_{s \in I(n,i)} |X_s'^{n,i} - X_{(i-1)/n}'^{n,i}|.$$

Standard arguments, using (3.9) and (3.10) and also the boundedness of $f$ for $\widetilde{X}_i^n$ and $\widetilde{X}_i'^n$, yield

$$\begin{aligned}
&E((\widehat{M}_i^n)^2|\mathcal{F}_{(i-1)/n}) \le C\frac{c_n}{n},\\
&E((\widehat{A}_i^n)^2|\mathcal{F}_{(i-1)/n}) \le C\frac{c_n^2 + d_n^2}{n},\\
&E((\widehat{Y}_i^n)^2|\mathcal{F}_{(i-1)/n}) \le C\left(\frac{c_n}{n} + \frac{d_n^2}{n^2}\right),\\
&E((\widetilde{X}_i^n)^2|\mathcal{F}_{(i-1)/n}) \le \frac{C}{n},\\
&E((\widetilde{X}_i'^n)^2|\mathcal{F}_{(i-1)/n}) \le C\left(\frac{c_n}{n} + \frac{d_n^2}{n^2}\right).
\end{aligned} \tag{3.16}$$

Let denote by $T(n,i)_p$ for $p = 1, 2, \dots,$ the successive jump times of $Y$, after $\frac{i-1}{n}$ and of size bigger than or equal to $\beta_n$. Let also $K(n,i)$ be the (random) integer such that $T(n,i)_{K(n,i)} \le \frac{i}{n} < T(n,i)_{K(n,i)+1}$. Then we set for $t \ge (i-1)/n$,

$$\begin{aligned}
V_t^{n,i} &= G(X_{(i-1)/n}, x)\mathbb{1}_{\{|x|>\beta_n\}}\mathbb{1}_{((i-1)/n,\infty)}(s) \star \mu_{T(n,i)_1 \wedge t},\\
V_t'^{n,i} &= N_t'^{\beta_n} - N_{(i-1)/n}'^{\beta_n}.
\end{aligned}$$

Now, we can introduce the variables $\zeta_i^n(j)$ occuring in (3.15):

$$\begin{aligned}
\zeta_i^n(1) &= u_n \int_{I(n,i)} (Y_{s-}'^{n,i} - g(X_{(i-1)/n})Y_{s-}^{n,i})\, dY_s,\\
\zeta_i^n(2) &= u_n g(X_{(i-1)/n})\left(\int_{I(n,i)} Y_{s-}^{n,i}\, dM_s^{\beta_n} + \int_{I(n,i)} M_{s-}^{n,i}\, dA_s^{\beta_n}\right),\\
\zeta_i^n(3) &= u_n \int_{I(n,i)} (V_{s-}'^{n,i} - V_{s-}^{n,i})\, dY_s^{\beta_n},\\
\zeta_i^n(4) &= u_n \int_{[T(n,i)_3, i/n]} V_{s-}'^{n,i}\, dN_s^{\beta_n} \mathbb{1}_{\{K(n,i)\ge 3\}},\\
\zeta_i^n(5) &= u_n(g(X_{T(n,i)_1-}) - g(X_{(i-1)/n}))\Delta Y_{T(n,i)_1} \Delta Y_{T(n,i)_2} \mathbb{1}_{\{K(n,i)\ge 2\}},\\
\zeta_i^n(6) &= u_n k(X_{T(n,i)_1-}, \Delta Y_{T(n,i)_1})\Delta Y_{T(n,i)_1}^2 \Delta Y_{T(n,i)_2} \mathbb{1}_{\{K(n,i)\ge 2\}},\\
\zeta_i^n(7) &= u_n g(X_{(i-1)/n})\left(\int_{I(n,i)} Y_{s-}^{n,i}\, dN_s^{\beta_n} - \Delta Y_{T(n,i)_1} Y_{T(n,i)_1}^{n,i} \mathbb{1}_{\{K(n,i)\ge 1\}}\right),
\end{aligned}$$

$$\zeta_i^n(8) = u_n k(X_{(i-1)/n}, \Delta Y_{T(n,i)_1})(\Delta Y_{T(n,i)_1})^2 (M^{n,i}_{i/n} - M^{n,i}_{T(n,i)_1}),$$

$$\zeta_i^n(9) = u_n g(X_{(i-1)/n}) \Delta Y_{T(n,i)_1} M^{n,i}_{i/n} \mathbb{1}_{\{K(n,i)\geq 1\}},$$

$$\zeta_i^n(10) = u_n g(X_{(i-1)/n}) \int_{I(n,i)} (A^{n,i}_{s-} + \Delta Y_{T(n,i)_1} \mathbb{1}_{\{K(n,i)\geq 1\}})\, dA^{\beta_n}_s,$$

$$\zeta_i^n(11) = u_n g(X_{(i-1)/n}) \Delta Y_{T(n,i)_1} \Delta Y_{T(n,i)_2} \mathbb{1}_{\{K(n,i)\geq 2\}},$$

$$\zeta_i^n(12) = u_n k(X_{(i-1)/n}, \Delta Y_{T(n,i)_1})(\Delta Y_{T(n,i)_1})^2\, d_n \left(\frac{i}{n} - T(n,i)_1\right) \mathbb{1}_{\{K(n,i)\geq 1\}},$$

$$\zeta_i^n(13) = M^{n,i}_{i/n} + \sum_{j\geq 2} \Delta Y_{T(n,i)_j} \mathbb{1}_{\{K(n,i)\geq j\}},$$

$$\zeta_i^n(14) = \frac{d_n}{n} + \Delta Y_{T(n,i)_1} \mathbb{1}_{\{K(n,i)\geq 1\}}.$$

Then we deduce from (3.7) and from (2.1) (after some tedious calculations) that (3.15) holds.

Finally, the following property will be used over and over again:

(3.17) Conditionally on $\mathcal{F}_{(i-1)/n}$ the variables $(\Delta Y_{T(n,i)_j})_{j\geq 1}, K(n,i), Y^{n,i}$ are independent; each $Y_{T(n,i)}$ has the law $\frac{1}{\theta(\beta_n)} F(dx) \mathbb{1}_{\{|x|>\beta_n\}}$; $K(n,i)$ has a Poisson law with parameter $\lambda_n$.

**4. The key lemma.** This section is devoted to proving the next lemma:

LEMMA 4.1. *Assume that we are in one of the Cases* 1, 2a, 2b, 3a *or* 3b.

(a) *We have* $\Gamma^n(j) \xrightarrow{P} 0$ *if* $j = 1, \dots, 8$, *and also*

- *for* $j = 9$ *in Cases* 2a, 2b, 3a *and* 3b;
- *for* $j = 10$ *in Cases* 1, 2b *and* 3b;
- *for* $j = 11$ *in Cases* 1 *and* 3a;
- *for* $j = 12$ *in Cases* 1, 2a, 2b *and* 3b;
- *for* $j = 13$ *in Cases* 1, 2b, 3a *and* 3b.

(b) *In the remaining cases, and for* $j = 9$ (*resp.* $j = 10$, *resp.* $j = 11$, *resp.* $j = 12$, *resp.* $j = 13$, *resp.* $j = 14$), *the sequences* $(\zeta_i^n(j))$ *satisfy* (2.9) [*resp.* (2.8), *resp.* (2.9), *resp.* (2.7), *resp.* (2.8), *resp.* (2.8)], *with* (2.11).

We will proceed through a large number of (very technical) steps.

4.1. *Step* 1: *auxiliary results.* Let us first derive some easy consequences of Cauchy–Schwarz and Doob inequalities. We consider a locally square-integrable martingale $N$ such that $\langle N, N\rangle_t = \int_0^t c_s\, ds$, where $c$ is a bounded process, and a constant $\gamma$, and also a bounded predictable process $H$, and fix $n$ and $i$, and set for $t > \frac{i-1}{n}$,

$$Z_t = \int_{(i-1)/n}^t H_s(\gamma\, ds + dN_s).$$

LEMMA 4.2. *In the above setting, we have for* $t > \frac{i-1}{n}$,

$$E(Z_t|\mathcal{F}_{(i-1)/n}) = \gamma \int_{(i-1)/n}^t E(H_s|\mathcal{F}_{(i-1)/n})\, ds, \tag{4.1}$$

$$\begin{aligned} E\Big(\sup_{t\in I(n,i)} Z_t^2|\mathcal{F}_{(i-1)/n}\Big) &\le \frac{2\gamma^2}{n^2} E\Big(\sup_{t\in I(n,i)} H_t^2\Big|\mathcal{F}_{(i-1)/n}\Big) \\ &\quad + \frac{8}{n} E\Big(\sup_{t\in I(n,i)} H_t^2 c_t\Big|\mathcal{F}_{(i-1)/n}\Big). \end{aligned} \tag{4.2}$$

Next we consider integrals w.r.t. the random measure $\mu$. Let $W$ be a predictable function on $\Omega \times \mathbb{R}_+ \times \mathbb{R}$, which is bounded on $\Omega \times \mathbb{R}_+ \times [-p, p]$ (recall that $p$ is such that $F$ charges only $[-p, p]$), and such that the process

$$H_s = \int_{\mathbb{R}} |W(s, x)| F(dx)$$

is also bounded. We fix again $n$ and $i$, and set for $t > \frac{i-1}{n}$,

$$Z_t = \int_{(i-1)/n}^t \int_{\mathbb{R}} W(s, x)\mu(ds, dx).$$

LEMMA 4.3. *In the above setting, we have for* $t > \frac{i-1}{n}$,

$$E(Z_t|\mathcal{F}_{(i-1)/n}) = \int_{(i-1)/n}^t ds\, E\Big(\int_{\mathbb{R}} W(s, x)F(dx)\Big|\mathcal{F}_{(i-1)/n}\Big), \tag{4.3}$$

$$\begin{aligned} E\Big(\sup_{t\in I(n,i)} Z_t^2\Big|\mathcal{F}_{(i-1)/n}\Big) &\le \frac{2}{n^2} E\Big(\sup_{t\in I(n,i)} \Big(\int_{\mathbb{R}} W(t, x)F(dx)\Big)^2\Big|\mathcal{F}_{(i-1)/n}\Big) \\ &\quad + \frac{8}{n} E\Big(\sup_{t\in I(n,i)} \int_{\mathbb{R}} W(t, x)^2 F(dx)\Big|\mathcal{F}_{(i-1)/n}\Big). \end{aligned} \tag{4.4}$$

PROOF. We can write $Z = Z' + Z''$, where

$$Z'_t = \int_{(i-1)/n}^t \int_{\mathbb{R}} W(s, x)\nu(ds, dx), \qquad Z''_t = Z_t - Z'_t.$$

Since $\nu(ds,dx)=ds\otimes F(dx)$ is the predictable compensator of $\mu$, we have $E(Z_t|\mathcal{F}_{(i-1)/n})=E(Z'_t|\mathcal{F}_{(i-1)/n})$, and (4.3) readily follows. We also get

$$E\Big(\sup_{t\in I(n,i)} Z'^2_t\Big|\mathcal{F}_{(i-1)/n}\Big)\le\frac{1}{n^2}E\Big(\sup_{t\in I(n,i)}\Big(\int_{\mathbb{R}}W(t,x)F(dx)\Big)^2\Big|\mathcal{F}_{(i-1)/n}\Big).$$

On the other hand, $Z''$ is a square-integrable martingale with bracket $\langle Z'',Z''\rangle_t=\int_{(i-1)/n}^t ds\int_{\mathbb{R}}W(s,x)^2F(dx)$ and, thus,

$$E\Big(\sup_{t\in I(n,i)} Z''^2_t\Big|\mathcal{F}_{(i-1)/n}\Big)\le 4\int_{(i-1)/n}^{i/n}ds\,E\Big(\int_{\mathbb{R}}W(t,x)^2F(dx)\Big|\mathcal{F}_{(i-1)/n}\Big).$$

Then (4.4) readily follows. $\square$

In the next auxiliary result, we fix $n$ and $i$ and write for simplicity $T_j=T(n,i)_j$ and $K=K(n,i)$. For some $r\ge 2$, let also $H_r$ be a random variables satisfying

$$\text{(4.5)}\qquad \begin{gathered} H_r,\ K \text{ and } \Delta Y_{T_r} \text{ are independent conditionally on } \mathcal{F}_{(i-1)/n},\\ |H_r|\le C_0\sum_{j=1}^{r-1}|\Delta Y_{T_j}|.\end{gathered}$$

Recall the notation (3.12) for $\lambda_n$. We also set $v_n(r)=e^{-\lambda_n}\frac{\lambda_n^r}{r!}$ and $H'_r=H_r\Delta Y_{T_r}$. The lemma below does not require any particular choice for the sequence $\beta_n$.

Lemma 4.4. *Under Hypothesis* (H1-$\alpha$) *with* $\alpha\in(0,2)$ *and* (4.5), *for all* $r'\ge r\ge 2$ *and* $y>0$, *we have the following estimates, where the constant* $C$ *below depends on* $C_0$ *in* (4.5):

$$\text{(4.6)}\qquad P(|u_nH'_r|>y,K=r'|\mathcal{F}_{(i-1)/n})\le C\frac{r^{1+\alpha}u_n^\alpha v_n(r')}{y^\alpha\theta(\beta_n)^2}\log\Big(\frac{1}{\beta_n}\Big),$$

$$\text{(4.7)}\qquad E(|u_nH_r|^\alpha\mathbb{1}_{\{K=r'\}}|\mathcal{F}_{(i-1)/n})\le Cr^2\frac{u_n^\alpha v_n(r')}{\theta(\beta_n)}\log\Big(\frac{1}{\beta_n}\Big),$$

$$\text{(4.8)}\qquad E(|u_nH'_r|^2\mathbb{1}_{\{|u_nH'_r|\le y,K=r'\}}|\mathcal{F}_{(i-1)/n})\le Cr^2y^{2-\alpha}\frac{u_n^\alpha v_n(r')}{\theta(\beta_n)^2}\log\Big(\frac{1}{\beta_n}\Big),$$

$$\text{(4.9)}\qquad E(|u_nH'_r|\mathbb{1}_{\{K=r'\}}|\mathcal{F}_{(i-1)/n})\le Cr\frac{u_nv_n(r')}{\theta(\beta_n)^2}\delta_n^2,$$

$$\text{(4.10)}\qquad \begin{aligned}\alpha>1\quad\Longrightarrow\quad &E(|u_nH'_r|\mathbb{1}_{\{|u_nH'_r|>y,K=r'\}}|\mathcal{F}_{(i-1)/n})\\ &\le C\frac{r^2u_n^\alpha v_n(r')}{y^{\alpha-1}\theta(\beta_n)^2}\log\Big(\frac{1}{\beta_n}\Big),\end{aligned}$$

$$\text{(4.11)} \qquad \alpha < 1 \quad \Longrightarrow \quad E(|u_n H'_r| \mathbb{1}_{\{|u_n H'_r| \le y, K = r'\}} | \mathcal{F}_{(i-1)/n}) \\ \le C \frac{r^2 u_n^\alpha v_n(r')}{y^{\alpha - 1} \theta(\beta_n)^2} \left( \log \left( \frac{1}{\beta_n} \right) \right)^2.$$

PROOF. Recalling (3.17) and (4.5), we see that the left-hand side of (4.6) is smaller than

$$\begin{aligned} &\le v_n(r') P\biggl( |\Delta Y_{T_r}| > \frac{y}{u_n C_0 \sum_{j=1}^{r-1} |\Delta Y_{T_j}|} \bigg| \mathcal{F}_{(i-1)/n} \biggr) \\ &= \frac{v_n(r')}{\theta(\beta_n)} \int_{\{|x| > \beta_n\}} F(dx) P\biggl( \sum_{j=1}^{r-1} |\Delta Y_{T_j}| > \frac{y}{u_n C_0 |x|} \bigg| \mathcal{F}_{(i-1)/n} \biggr) \\ &\le \frac{r v_n(r')}{\theta(\beta_n)} \int_{\{|x| > \beta_n\}} F(dx) P\biggl( |\Delta Y_{T_1}| > \frac{y}{(r-1) u_n C_0 |x|} \bigg| \mathcal{F}_{(i-1)/n} \biggr) \\ &\le \frac{r v_n(r')}{\theta(\beta_n)^2} \int_{\{|x| > \beta_n\}} F(dx) E\biggl( \int F(dz) \mathbb{1}_{\{|z| > y/(r-1) u_n C_0 |x|\}} \bigg| \mathcal{F}_{(i-1)/n} \biggr) \\ &\le C \frac{r^{1+\alpha} u_n^\alpha v_n(r')}{\theta(\beta_n)^2 y^\alpha} \rho_n \le C \frac{r^{1+\alpha} u_n^\alpha v_n(r')}{y^\alpha \theta(\beta_n)^2} \log\left( \frac{1}{\beta_n} \right), \end{aligned}$$

where we have used Hypothesis (H1-$\alpha$) and (3.3). Hence, (4.6) holds.

Second,

$$\begin{aligned} E(|H_r|^\alpha \mathbb{1}_{\{K = r'\}} | \mathcal{F}_{(i-1)/n}) &\le C v_n(r') (r-1)^{\alpha \vee 1} E(|\Delta Y_{T_1}|^\alpha | \mathcal{F}_{(i-1)/n}) \\ &\le r^2 \frac{v_n(r')}{\theta(\beta_n)} \rho_n, \end{aligned}$$

and (4.7) follows from (3.3).

Third, by (3.3) the left-hand side of (4.8) is

$$\begin{aligned} &= E((u_n H_r)^2 |\Delta Y_{T_r}|^2 \mathbb{1}_{\{|\Delta Y_{T_r}| \le y / |u_n H_r|\}} \mathbb{1}_{\{K = r'\}} | \mathcal{F}_{(i-1)/n}) \\ &= \frac{1}{\theta(\beta_n)} E\biggl( (u_n H_r)^2 c\biggl( \frac{y}{|u_n H_r|} \biggr) \mathbb{1}_{\{K = r'\}} \bigg| \mathcal{F}_{(i-1)/n} \biggr) \\ &\le C \frac{y^{2-\alpha}}{\theta(\beta_n)} E(|u_n H_r|^\alpha \mathbb{1}_{\{K = r'\}} | \mathcal{F}_{(i-1)/n}), \end{aligned}$$

and (4.8) follows from (4.7).

Next the left-hand side of (4.9) is smaller than $C u_n (r-1) v_n(r') (E(|\Delta Y_{T_1}|))^2$. But $E(|\Delta Y_{T_1}|) = \frac{\delta_n}{\theta(\beta_n)}$, hence, (4.9).

Next, suppose that $\alpha > 1$. The left-hand side of (4.10) is

$$= E(|u_n H_r| |\Delta Y_{T_r}| \mathbb{1}_{\{|\Delta Y_{T_r}| > y / |u_n H_r|\}} \mathbb{1}_{\{K = r'\}} | \mathcal{F}_{(i-1)/n})$$

$$\leq \frac{1}{\theta(\beta_n)} E\biggl(|u_n H_r| \delta\biggl(\frac{y}{|u_n H_r|}\biggr) \mathbb{1}_{\{K=r'\}} \bigg| \mathcal{F}_{(i-1)/n}\biggr)$$

$$\leq \frac{C}{y^{\alpha-1}\theta(\beta_n)} E(|u_n H_r|^{\alpha} \mathbb{1}_{\{K=r'\}} | \mathcal{F}_{(i-1)/n})$$

[use (3.3)], and (4.10) follows from (4.7).

Finally, suppose that $\alpha < 1$. The left-hand side of (4.11) is

$$= E(|u_n H_r||\Delta Y_{T_r}| \mathbb{1}_{\{|\Delta Y_{T_r}| \leq y/|u_n H_r|\}} \mathbb{1}_{\{K=r'\}} | \mathcal{F}_{(i-1)/n})$$

$$\leq E\biggl(|u_n H_r||\Delta Y_{T_r}|^{\alpha} \frac{y^{1-\alpha}}{|u_n H_r|^{1-\alpha}} \mathbb{1}_{\{K=r'\}} \bigg| \mathcal{F}_{(i-1)/n}\biggr)$$

$$\leq \frac{C\rho_n y^{1-\alpha}}{\theta(\beta_n)} E(|u_n H_r|^{\alpha} \mathbb{1}_{\{K=r'\}} | \mathcal{F}_{(i-1)/n})$$

and (4.11) follows from (3.3) and (4.7). $\square$

Now we proceed to proving Lemma 4.1, going step by step.

4.2. *Step* 2: $j = 1, 2, 3$. First, from (3.8) and (3.11) we see that for $t > \frac{i-1}{n}$,

$$Y_t'^{n,i} - g(X_{(i-1)/n}) Y_t^{n,i}$$
$$= \int_{(i-1)/n}^{t} H_s(d_n\, ds + dM_s^{\beta_n}) + \int_{(i-1)/n}^{t} \int_{\mathbb{R}} W(s,x)\mu(ds,dx),$$

where $H_t = g(X_{t-}) - g(X_{(i-1)/n})$ and $W(t,x) = x^2 k(X_{t-}, x) \mathbb{1}_{\{|x| \leq \beta_n\}}$ (with $n$ and $i$ fixed). Then since $g$ is bounded with a bounded derivative and $k$ is bounded over $\mathbb{R} \times [-p, p]$, we deduce from (4.2) (with $N = M^{\beta_n}$) and (4.4), together with (3.16) and the fact that $|H_t| \leq C\widetilde{X}_i^n$ for $t \in I(n,i)$, that

$$E\Bigl(\sup_{s \in I(n,i)} (Y_s'^{n,i} - g(X_{(i-1)/n}) Y_s^{n,i})^2 \Big| \mathcal{F}_{(i-1)/n}\Bigr) \leq C\alpha_n$$
$$\text{where } \alpha_n = \frac{c_n \beta_n^2}{n} + \frac{c_n}{n^2} + \frac{d_n^2}{n^3}.$$

Now we have $Y = M^p + A^p$ and $A_t^p = b't$: we can apply once more Lemma 4.2 (with $N = M^p$) to get (2.8) for $\zeta_i^n(1)$ with $\xi_n = Cu_n\sqrt{\alpha_n}$ and $\xi_n' = C\xi_n^2$. Since $\xi_n \to 0$ by (3.13), we obtain the result for $j = 1$ from Lemma 2.5.

By (3.16) and Lemma 4.2 again, we have (2.8) for $\zeta_i^n(2)$ with $\xi_n = 0$ and $\xi_n' = Cu_n^2 c_n(\frac{c_n}{n} + \frac{d_n^2}{n^2})$, which goes to 0 by (3.13): hence, the result for $j = 2$.

Next we have $V'^{n,i} - V^{n,i} = U^{n,i} \star \mu$ [use (3.8)], where

$$U^{n,i}(s,x) = (g(X_{s-})\mathbb{1}_{\{s>T(n,i)_1\}} + (g(X_{s-}) - g(X_{(i-1)/n}))\mathbb{1}_{\{(i-1)/n<s\le T(n,i)_1\}})x\mathbb{1}_{\{|x|>\beta_n\}} + (k(X_{s-},x)\mathbb{1}_{\{s>T(n,i)_1\}} + (k(X_{s-},x) - k(X_{(i-1)/n},x))\mathbb{1}_{\{(i-1)/n<s\le T(n,i)_1\}})x^2\mathbb{1}_{\{|x|>\beta_n\}}.$$

The fact that $f$ is $C^3$ with compact support implies that $|g(x)-g(x')| \le Cu$ and $|k(x,y)-k(x',y)| \le Cu$ whenever $|x-x'| \le u$ and $y \in [-p,p]$. Then if $t \in I(n,i)$, we have

$$\left|\int_{\mathbb{R}} U^{n,i}(t,x)F(dx)\right| \le C(|d'_n|+1)(\mathbb{1}_{\{K(n,i)\ge 1\}} + \widetilde{X}^n_i),$$

$$\int_{\mathbb{R}} U^{n,i}(t,x)^2 F(dx) \le C(\mathbb{1}_{\{K(n,i)\ge 1\}} + (\widetilde{X}^n_i)^2).$$

Then it follows from (3.16), (3.17) and from Lemma 4.3 that

$$\sup_{s\in I(n,i)} |E(V'^{n,i}_s - V^{n,i}_s|\mathcal{F}_{(i-1)/n})| \le C\alpha'_n$$

$$\text{where } \alpha'_n = \frac{1}{n}(1+|d_n|)\left(\lambda_n + \frac{1}{\sqrt{n}}\right),$$

$$E\left(\sup_{s\in I(n,i)} |V'^{n,i}_s - V^{n,i}_s|^2 \Big| \mathcal{F}_{(i-1)/n}\right) \le C\alpha''_n$$

$$\text{where } \alpha''_n = \frac{1}{n}\left(1+\frac{d_n^2}{n}\right)\left(\lambda_n + \frac{1}{n}\right).$$

Recall that $Y^{\beta_n} = M^{\beta_n} + A^{\beta_n}$ and $A^{\beta_n}_t = d_n t$, so the above estimates and an application of Lemma 4.2 (with $N = M^{\beta_n}$) allow us to deduce that $\zeta^n_i(3)$ satisfies (2.8) with $\xi_n = Cu_n|d_n|\alpha'_n$ and $\xi'_n = Cu_n^2\alpha''_n(c_n + d_n^2/n)$. By (3.13) these go to 0, hence, the result for $j=3$.

4.3. *Step* 3: $j=4$. In order to study $\zeta^n_i(4)$ we apply Lemma 4.4: we use the notation $T_j$ and $K$ of Step 1 and we set $H_r = \Delta Y'_{T_1} + \cdots + \Delta Y'_{T_{r-1}}$ and $H'_r = H_r \Delta Y_{T_r}$ and $H''_k = u_n \sum_{r=3}^k H'_r$ (an empty sum being set to 0). The key observation is then that

$$\zeta^n_i(4) = H''_K \mathbb{1}_{\{K\ge 3\}}.$$

Observe that each $H_r$ satisfies (4.5), and recall (3.17). We will also use the easily proven fact that for any $a, a' > 0$ and $r \in \mathbb{N}$, there is a constant $C_{a,a',r}$ such that (recall that $\lambda_n \to 0$)

$$0 < \lambda_n \le \frac{1}{a} \implies \sum_{k=r}^{\infty} v_n(k) a^k k^{a'} \le C_{a,a',r}\lambda_n^r. \tag{4.12}$$

In view of (4.6), we get for $y>0$ and $k'\ge k\ge 3$,

$$P(|H''_k|>y,K=k'|\mathcal{F}_{(i-1)/n})\le \sum_{r=3}^{k} P\biggl(|u_nH'_r|>\frac{y}{k-2},K=k'\Bigm|\mathcal{F}_{(i-1)/n}\biggr)$$

(4.13)

$$\le C\frac{k^{2+2\alpha}u_n^\alpha v_n(k')}{y^\alpha\theta(\beta_n)^2}\log\biggl(\frac{1}{\beta_n}\biggr).$$

Therefore, since $\lambda_n=\theta(\beta_n)/n$, we deduce from (4.12) and (4.13),

$$P(|\zeta_i^n(4)|>y|\mathcal{F}_{(i-1)/n})=\sum_{k=3}^{\infty}P(|H''_k|>y,K=k|\mathcal{F}_{(i-1)/n})$$

(4.14)

$$\le \frac{Cu_n^\alpha}{n^2y^\alpha}\log\biggl(\frac{1}{\beta_n}\biggr)\lambda_n.$$

Next, a simple computation shows that for $k\ge 3$,

$$H''^2_k\mathbb{1}_{\{|H''_k|\le 1\}}\le \sum_{r=3}^{k-1}2^{k-r}\mathbb{1}_{\{|H''_r|>1\}}+\sum_{r=3}^{k}2^{k+1-r}(u_nH'_r)^2\mathbb{1}_{\{|u_nH'_r|\le 2\}}.$$

Hence, we obtain, by virtue of (4.8) and (4.13),

$$E(H''^2_k\mathbb{1}_{\{|H''_k|\le 1,K=k\}}|\mathcal{F}_{(i-1)/n})\le C2^kk^{2\alpha+3}\frac{u_n^\alpha v_n(k)}{\theta(\beta_n)^2}\log\biggl(\frac{1}{\beta_n}\biggr),$$

and by (4.12),

$$E(\zeta_i^n(4)^2\mathbb{1}_{\{|\zeta_i^n(4)|\le 1\}}|\mathcal{F}_{(i-1)/n})=\sum_{k=3}^{\infty}E(|H''_k|^2\mathbb{1}_{\{|H''_k|\le 1,K=k\}}|\mathcal{F}_{(i-1)/n})$$

(4.15)

$$\le \frac{Cu_n^\alpha}{n^2}\log\biggl(\frac{1}{\beta_n}\biggr)\lambda_n.$$

Next, for $r\ge 2$, we have $H''_r\mathbb{1}_{\{|H''_r|\le y|\}}=\sum_{j=1}^4\mu_r^j(y)$, where

$$\mu_r^1(y)=H''_{r-1}\mathbb{1}_{\{|H''_{r-1}|\le y/2\}},\qquad \mu_r^2(y)=-H''_{r-1}\mathbb{1}_{\{|H''_{r-1}|\le y/2,|H''_r|>y\}},$$

$$\mu_r^3(y)=H''_r\mathbb{1}_{\{|H''_r|\le y,|H''_{r-1}|>y/2\}},\qquad \mu_r^4(y)=u_nH'_r\mathbb{1}_{\{|H''_{r-1}|\le y/2,|H''_r|\le y\}}.$$

Inequality (4.13) yields, for $k\ge r$,

$$E(|\mu_r^2(y)+\mu_r^3(y)|\mathbb{1}_{\{K=k\}}|\mathcal{F}_{(i-1)/n})\le Cr^{2+2\alpha}\frac{u_n^\alpha v_n(k)}{\theta(\beta_n)^2}y^{1-\alpha}\log\biggl(\frac{1}{\beta_n}\biggr).$$

(4.16)

Set also

$$\nu_n=\begin{cases}\biggl(\log\biggl(\dfrac{1}{\beta_n}\biggr)\biggr)^2, & \text{if } \alpha\le 1,\\ u_n^{1-\alpha}\delta_n^2, & \text{if } \alpha>1.\end{cases} \tag{4.17}$$

Note that $|\mu_r^4(y)| \le |u_n H'_r| \mathbb{1}_{\{|u_n H'_r| \le 3y/2\}}$, so (4.9) when $\alpha \ge 1$ and (4.11) when $\alpha < 1$ yield

$$(4.18)\quad k \ge r \quad \Longrightarrow \quad |E(\mu_r^4(y)\mathbb{1}_{\{K=k\}}|\mathcal{F}_{(i-1)/n})| \le Cr^2 \frac{u_n^\alpha \nu_n v_n(k)}{\theta(\beta_n)^2}(1+y^{1-\alpha}).$$

Putting (4.16) and (4.18) together, and setting $\xi_{r,k}(y) = |E(H''_r \mathbb{1}_{\{|H''_r| \le y, K=k\}}|\mathcal{F}_{(i-1)/n})|$, we get for $3 \le r \le k$,

$$\xi_{r,k}(y) \le \xi_{r-1,k}\left(\frac{y}{2}\right) + Cr^{2+2\alpha}(1+y^{1-\alpha})\frac{u_n^\alpha v_n(k)}{\theta(\beta_n)^2}\left(\nu_n + \log\left(\frac{1}{\beta_n}\right)\right).$$

Recalling that $H''_2 = 0$, hence, $\xi_{2,k}(y) = 0$, an induction gives

$$\xi_{r,k}(1) \le Cr^{3+2\alpha}2^{r(\alpha-1)^+}\frac{u_n^\alpha v_n(k)}{\theta(\beta_n)^2}\left(\nu_n + \log\left(\frac{1}{\beta_n}\right)\right)$$

and, thus, we obtain by (4.12),

$$(4.19)\quad \begin{aligned} |E(H''_K \mathbb{1}_{\{|H''_K| \le 1, K \ge 3\}}|\mathcal{F}_{(i-1)/n})| &\le \sum_{k=3}^{\infty} \xi_{k,k}(1) \\ &\le C\frac{u_n^\alpha}{n^2}\left(\nu_n + \log\left(\frac{1}{\beta_n}\right)\right)\lambda_n, \\ |E(H''_K \mathbb{1}_{\{|H''_K| \le 1, K \ge 4\}}|\mathcal{F}_{(i-1)/n})| &\le \sum_{k=4}^{\infty} \xi_{k,k}(1) \\ &\le C\frac{u_n^\alpha}{n^2}\left(\nu_n + \log\left(\frac{1}{\beta_n}\right)\right)\lambda_n^2. \end{aligned}$$

In Cases 2b and 3b, the symmetry property of $\Delta Y_{T_3}$ and (3.17) and the fact that $H''_3 = u_n H_3 \Delta Y_{T_3}$ yield that $E(H''_3 \mathbb{1}_{\{|H''_3| \le 1, K=3\}}|\mathcal{F}_{(i-1)/n}) = 0$; therefore, the above two inequalities yield, with $\lambda'_n = \lambda_n$ in Cases 2b and 3b and $\lambda'_n = 1$ otherwise,

$$(4.20)\quad |E(\zeta_i^n(4)\mathbb{1}_{\{|\zeta_i^n(4)| \le 1\}}|\mathcal{F}_{(i-1)/n})| \le C\frac{u_n^\alpha}{n^2}\left(\nu_n + \log\left(\frac{1}{\beta_n}\right)\right)\lambda_n \lambda'_n.$$

Then we put together (4.14), (4.15) and (4.20): we see that $\zeta_i^n(4)$ satisfies (2.9) with

$$\xi_n = \frac{Cu_n^\alpha}{n}\left(\nu_n + \log\left(\frac{1}{\beta_n}\right)\right)\lambda_n\lambda'_n,$$

$$\xi'_n = \frac{Cu_n^\alpha}{n}\log\left(\frac{1}{\beta_n}\right)\lambda_n, \qquad \xi''_{n,y} = \frac{Cu_n^\alpha}{ny^\alpha}\log\left(\frac{1}{\beta_n}\right)\lambda_n.$$

Using (3.13), we see that the above quantities go to 0 as $n \to \infty$, hence, the result for $j = 4$.

4.4. *Step* 4: $j=5$. Observe that $|\zeta_i^n(5)| \leq C u_n \widetilde{X}_i'^n |\Delta Y_{T(n,i)_1} \Delta Y_{T(n,i)_2}| \mathbb{1}_{\{K(n,i)\geq 2\}}$. Then by (3.16) and (3.17),

$$E(|\zeta_i^n(5)||\mathcal{F}_{(i-1)/n}) \leq C u_n \frac{P(K(n,i)\geq 2|\mathcal{F}_{(i-1)/n})}{\theta(\beta_n)^2}\left(\frac{c_n}{n}+\frac{d_n^2}{n^2}\right)^{1/2}\delta_n^2$$
$$\leq \frac{C}{n}\frac{u_n(c_n/n+d_n^2/n^2)^{1/2}\delta_n^2}{n},$$
$$E(\zeta_i^n(5)^2|\mathcal{F}_{(i-1)/n}) \leq C u_n^2 \frac{P(K(n,i)\geq 2|\mathcal{F}_{(i-1)/n})}{\theta(\beta_n)^2}\left(\frac{c_n}{n}+\frac{d_n^2}{n^2}\right)$$
$$\leq \frac{C}{n}\frac{u_n^2}{n}\left(\frac{c_n}{n}+\frac{d_n^2}{n^2}\right).$$

Hence, (2.8) holds with $\xi_n = C\frac{u_n\delta_n^2}{n}(\frac{c_n}{n}+\frac{d_n^2}{n^2})^{1/2}$ and $\xi_n' = C\frac{u_n^2}{n}(\frac{c_n}{n}+\frac{d_n^2}{n^2})$. These sequences go to 0 by (3.13), hence, the result for $j=5$.

4.5. *Step* 5: $j=6$. Set $U_i^n = u_n \Delta Y_{T(n,i)_1}^2 |\Delta Y_{T(n,i)_2}| \mathbb{1}_{\{K(n,i)\geq 2\}}$. As for the proof of (4.6), and using (3.17), we get,

$$\begin{aligned}
&P(U_i^n > y|\mathcal{F}_{(i-1)/n})\\
&\quad\leq C\lambda_n^2 P\left(\Delta Y_{T(n,i)_1}^2|\Delta Y_{T(n,i)_2}| > \frac{y}{u_n}\Big|\mathcal{F}_{(i-1)/n}\right)\\
&\quad\leq \frac{C}{n^2}\int_{\{|x|>\beta_n\}} F(dx)\int_{\{|x'|>y/u_n x^2\}} F(dx')\\
&\quad\leq \frac{C}{n^2}\int F(dx)\frac{u_n^\alpha |x|^{2\alpha}}{y^\alpha} \leq C\frac{u_n^\alpha}{n^2 y^\alpha}
\end{aligned}\tag{4.21}$$

because $\int |x|^{2\alpha}F(dx) < \infty$. Similarly, we have by (3.3),

$$\begin{aligned}
&P(U_i^n \mathbb{1}_{\{U_i^n\leq 1\}}\mathcal{F}_{(i-1)/n})\\
&\quad\leq \frac{Cu_n}{n^2}\int_{\{|x|>\beta_n\}} F(dx)|x|\int_{\{|x'|\leq 1/\sqrt{u_n|x|}\}} x'^2 F(dx')\\
&\quad\leq \frac{Cu_n^{\alpha/2}}{n^2}\int_{\{|x|>\beta_n\}} F(dx)|x|^{\alpha/2} \leq \frac{Cu_n^{\alpha/2}}{n^2\beta_n^{\alpha/2}}.
\end{aligned}\tag{4.22}$$

Now we observe that $|k(x,y)| \leq C_0$ for some constant $C_0>0$, hence, $|\zeta_i^n(6)| \leq C_0 U_i^n$. Then if $|\zeta_i^n(6)| > y$ for some $y>0$, we must have $U_i^n > y/C_0$; also if $|\zeta_i^n(6)|\leq 1$, then we have $\zeta_i^n(6)^2 \leq |\zeta_i^n(6)| \leq C_0 U_i^n \mathbb{1}_{\{U_i^n\leq 1\}} + \mathbb{1}_{\{U_i^n>1\}}$. Then it readily follows from (4.21) and (4.22) that the sequence $(\zeta_i^n(6))$ satisfies (2.9) with $\xi_n = \xi_n' = \frac{C}{n}(u_n^\alpha + \frac{u_n^{\alpha/2}}{\beta_n^{\alpha/2}})$ and $\xi_{n,y}'' = \frac{Cu_n^\alpha}{ny^\alpha}$. Those sequences all go to 0, hence, the result for $j=6$.

4.6. *Step* 6: $j=7$. We use again all the notation of Step 1, so that (since $Y^{n,i}$ does not jump at times $T_j$)

$$\zeta_i^n(7)=u_n g(X_{(i-1)/n})\sum_{k=2}^K H_k, \qquad \text{where } H_k=Y^{n,i}_{T_k\wedge i/n}\Delta Y_{T_k}.$$

On the one hand we have $H_k^2\le(\widehat{Y}_i^n)^2\Delta Y_{T_k}^2$, so

$$\begin{aligned} E(\zeta_i^n(7)^2|\mathcal{F}_{(i-1)/n}) &\le Cu_n^2 E\Bigg(K\sum_{k=2}^K \Delta Y_{T_k}^2(\widehat{Y}_i^n)^2\Big|\mathcal{F}_{(i-1)/n}\Bigg)\\ &\le Cu_n^2\delta_n\sum_{k=2}^\infty \frac{v_n(k)}{\theta(\beta_n)}\int_{\{|x|>\beta_n\}}x^2F(dx) \qquad (4.23)\\ &\le C\frac{u_n^2}{n}\lambda_n\bigg(\frac{c_n}{n}+\frac{d_n^2}{n^2}\bigg)\end{aligned}$$

because of (3.17), (3.16) and (4.12). On the other hand we can write

$$\begin{aligned} &E(\zeta_i^n(7)|\mathcal{F}_{(i-1)/n})\\ &\quad =u_n g(X_{(i-1)/n})\sum_{k=2}^\infty E(\Delta Y_{T_k}\mathbb{1}_{\{T_k\le i/n\}}Y^{n,i}_{T_k}|\mathcal{F}_{(i-1)/n})\\ &\quad =u_n g(X_{(i-1)/n})\,d_n\sum_{k=2}^\infty E\bigg(\Delta Y_{T_k}\mathbb{1}_{\{T_k\le i/n\}}\bigg(T_k-\frac{i-1}{n}\bigg)\bigg|\mathcal{F}_{(i-1)/n}\bigg),\end{aligned}$$

again by (3.17) and because $Y^{n,i}$ is equal to a martingale plus $d_n(t-\frac{i-1}{n})$. Therefore,

$$\begin{aligned} |E(\zeta_i^n(7)|\mathcal{F}_{(i-1)/n})| &\le C\frac{u_n|d_n|}{n}E\Bigg(\sum_{k=2}^K|\Delta Y_{T_k}|\Big|\mathcal{F}_{(i-1)/n}\Bigg)\\ &= C\frac{u_n|d_n|}{n}\sum_{k=2}^\infty\frac{v_n(k)}{\theta(\beta_n)}\int_{\{|x|>\beta_n\}}|x|F(dx) \qquad (4.24)\\ &\le C\frac{u_n|d_n|}{n^2}\delta_n\lambda_n.\end{aligned}$$

So $\zeta_i^n(7)$ satisfies (2.8) with $\xi_n=C\frac{u_n|d_n|\delta_n\lambda_n}{n}$ and $\xi'_n=Cu_n^2(\frac{c_n}{n}+\frac{d_n^2}{n^2})\lambda_n$. Those sequences go to 0, hence, the result for $j=7$.

4.7. *Step* 7: $j=8,9$. We have $E(\zeta_i^n(8)|\mathcal{F}_{(i-1)/n})=E(\zeta_i^n(9)|\mathcal{F}_{(i-1)/n})=0$ [use (3.17)], and also $|\zeta_i^n(8)|+|\zeta_i^n(9)|\le Cu_n|\Delta Y_{T(n,i)_1}|\widehat{M}_i^n\mathbb{1}_{\{K(n,i)\ge1\}}$, hence, by (3.17) again and (3.16),

$$E(\zeta_i^n(8)^2|\mathcal{F}_{(i-1)/n})+E(\zeta_i^n(9)^2|\mathcal{F}_{(i-1)/n})\le\frac{Cu_n^2c_n}{n\theta(\beta_n)}P(K\ge1)\le\frac{Cu_n^2c_n}{n^2}.$$

Then (2.8) holds for $j=8$ and $j=9$, with $\xi_n=0$ and $\xi'_n=Cu_n^2c_n/n$. Except in Case 1, we have $\xi'_n\to 0$: hence, the result for $j=8$ and $j=9$ if we are not in Case 1.

It remains to study Case 1. For this we will use Lemma 2.6. We first calculate the conditional characteristic function $\phi_{n,i}(v)=E(e^{iv\zeta_i^n(8)}|\mathcal{F}_{(i-1)/n})$ and $\phi'_{n,i}(v)=E(e^{iv\zeta_i^n(9)}|\mathcal{F}_{(i-1)/n})$. Recall that $M^{n,i}$ is a Lévy process, independent of $\mathcal{F}_{(i-1)/n}$ and satisfying

$$E(e^{ivM_t^{n,i}})=\exp t\int_{\{|x|\le\beta_n\}}(e^{ivy}-1-ivy)F(dy). \tag{4.25}$$

Then, using (3.17) and the form of $\zeta_i^n(8)$ and $\zeta_i^n(9)$, we see that

$$\begin{aligned}\phi_{n,i}(v)&=e^{-\lambda_n}\\&\qquad+\int_0^{1/n}e^{-\theta(\beta_n)s}\,ds\int_{\{|x|>\beta_n\}}F(dx)E(e^{ivu_nx^2k(X_{(i-1)/n},x)(M_{1/n}^{n,i}-M_s^{n,i})}),\\ \phi'_{n,i}(v)&=e^{-\lambda_n}+\frac{1-e^{-\lambda_n}}{n\lambda_n}\int_{\{|x|>\beta_n\}}F(dx)E(e^{ivu_ng(X_{(i-1)/n})xM_s^{n,i}}).\end{aligned}$$

Observe that $1-e^{-\lambda_n}=\int_0^{1/n}e^{-\theta(\beta_n)s}\,ds\int_{\{|x|>\beta_n\}}F(dx)$. Then if

$$\begin{aligned}\gamma_{n,i}(v,x,y)&=e^{ivu_nx^2k(X_{(i-1)/n},x)y}-1-ivu_nx^2k(X_{(i-1)/n},x)y,\\ z_{n,i}(x,v,t)&=\frac{1-t}{n}\int_{\{|y|\le\beta_n\}}\gamma_{n,i}(v,x,y)F(dy),\end{aligned}$$

the change of variable $t=ns$ gives

$$\phi_{n,i}(v)=1+\frac{1}{n}\int_0^1e^{-\lambda_nt}\,dt\int_{\{|x|>\beta_n\}}F(dx)(e^{z_{n,i}(x,v,t)}-1). \tag{4.26}$$

In a similar way we obtain

$$\phi'_{n,i}(v)=1+\frac{1-e^{-\lambda_n}}{n\lambda_n}\int_{\{|x|>\beta_n\}}F(dx)(e^{z'_{n,i}(x,v)}-1), \tag{4.27}$$

where

$$\begin{aligned}\gamma'_{n,i}(v,x,y)&=e^{ivu_ng(X_{(i-1)/n})xy}-1-ivu_ng(X_{(i-1)/n})xy,\\ z_{n,i}(x,v)&=\frac{1}{n}\int_{\{|y|\le\beta_n\}}\gamma'_{n,i}(v,x,y)F(dy).\end{aligned}$$

Since $g$ (resp. $k$) is bounded (resp. bounded on $\mathbb{R}\times[-p,p]$), we get $|\gamma_{n,i}(v,x,y)|\le C(|vu_nx^2y|\wedge|vu_nx^2y|^2)$ and $|\gamma'_{n,i}(v,x,y)|\le C(|vu_nx^2y|\wedge|vu_nx^2y|^2)$

whenever $|x| \le p$. But (3.3) yields

$$\int (|uy| \wedge |uy|^2) F(dy) = |u| \int_{\{|y|>1/u\}} |y| F(dy) + u^2 \int_{\{|y|\le 1/u\}} y^2 F(dy) \le C|u|^\alpha.$$

Therefore, $|z_{n,i}(x,v,t)| \le \frac{C}{n}|vu_n x^2|^\alpha$ and $|z'_{n,i}(x,v)| \le \frac{C}{n}|vu_n x|^\alpha$, provided $t \in [0,1]$ and $|x| \le p$. In particular, the suprema of $|z_{n,i}|$ and $|z'_{n,i}|$ over all $t \in [0,1]$, $i$, $x \in [-p,p]$ and $v \in [-1,1]$ go to 0 as $n \to \infty$ and are thus uniformly bounded in $n$ as well. Therefore, in (4.26) and (4.27) the term $e^{\cdot} - 1$ is smaller than $\frac{C}{n}|vu_n x^2|^\alpha$ and $\frac{C}{n}|vu_n x|^\alpha$, respectively. Since $x \mapsto |x|^{2\alpha}$ is $F$-integrable, and using also (3.3), we readily deduce

$$|\phi_{n,i}(v) - 1| \le \frac{C|v|^\alpha u_n^\alpha}{n^2}, \qquad |\phi'_{n,i}(v) - 1| \le \frac{C|v|^\alpha u_n^\alpha}{n^2} \log\frac{1}{\beta_n}.$$

In other words, the sequences $(\zeta_i^n(8))$ and $(\zeta_i^n(9))$ satisfy (2.12) with, respectively, $\xi'''_{n,v} = Cv^\alpha u_n^\alpha/n$ and $\xi'''_{n,v} = Cv^\alpha u_n^\alpha(\log\frac{1}{\beta_n})/n$. In the first case we have $\xi'''_{n,v} \to 0$ for all $v \le 1$; in the second case we have $\xi'''_{n,v} \le Cv^\alpha$: then, combining Lemmas 2.6 and 2.5, we obtain the result for $j = 8$ and $j = 9$ in Case 1.

4.8. *Step* 8: $j = 10, 11, 12$. First consider $j = 10$. We have

$$\zeta_i^n(10) = u_n g(X_{(i-1)/n}) \frac{d_n}{n}\left(\frac{d_n}{2n} + \Delta Y_{T(n,i)_1} \mathbb{1}_{\{K(n,i)\ge 1\}}\right).$$

Then $a_i^n = E(\zeta_i^n(10)|\mathcal{F}_{(i-1)/n})$ and $b_i^n = E(\zeta_i^n(10)^2|\mathcal{F}_{(i-1)/n})$ have

$$a_i^n = u_n g(X_{(i-1)/n}) \frac{d_n}{n}\left(\frac{d_n}{2n} + (1 - e^{-\lambda_n}) \frac{d'_n}{\theta(\beta_n)}\right), \tag{4.28}$$

$$\begin{aligned} b_i^n &\le \frac{C u_n^2 d_n^2}{n^2}\left(\frac{d_n^2}{n^2} + E(\Delta Y^2_{T(n,i)_1} \mathbb{1}_{\{K(n,i)\ge 1\}}|\mathcal{F}_{(i-1)/n})\right) \\ &\le C\frac{u_n^2 d_n^2}{n^3}\left(1 + \frac{d_n^2}{n}\right). \end{aligned} \tag{4.29}$$

In particular, $|a_i^n| \le C\frac{u_n(d_n^2 + |d_n d'_n|)}{n^2}$, so (2.8) holds with $\xi_n = C\frac{u_n(d_n^2 + |d_n d'_n|)}{n}$ and $\xi'_n = C\frac{u_n^2 d_n^2}{n^2}(1 + \frac{d_n^2}{n})$: these quantities go to 0 in Cases 1, 2b and 3a, and they are always bounded: hence, the result for $j = 10$.

Next, consider $j = 11$. We apply Lemma 4.4, in which we set $H_2 = g(X_{(i-1)/n})\Delta Y_{T_1}$, so that (4.5) is satisfied and $\zeta_i^n(11) = u_n H'_2 \mathbb{1}_{\{K\ge 2\}}$. Then using (4.6) and (4.8), and also (4.11) [resp. (4.9)] if $\alpha < 1$ (resp. if $\alpha \ge 1$) for $H_2'$ and summing up over $r'\geq 2$, we get that when $\alpha\leq 1$, the sequences $(\zeta_i^n(11))$ satisfy (2.9) with

$$\xi_n = C\frac{u_n^\alpha}{n}\left(\log\frac{1}{\beta_n}\right)^2, \qquad \xi_n' = C\frac{u_n^\alpha}{n}\log\frac{1}{\beta_n}, \qquad \xi_{n,y}'' = C\frac{u_n^\alpha}{ny^\alpha}\log\frac{1}{\beta_n}. \tag{4.30}$$

In Case 3a (resp. 2a) we have $\xi_n'\to 0$ and $\xi_{n,y}''\to 0$ and also $\xi_n\to 0$ (resp. $\xi_n\leq C$): hence, the result for $j=11$ in Cases 2a and 3a. In Cases 2b and 3b we still have (4.30) but, for symmetry reasons, we may take $\xi_n=0$: hence, we also get the result in Cases 2b and 3b.

Now consider Case 1. The estimates (4.30) are not fine enough and we have to resort on another method. The $\zeta_i^n(11)$'s satisfy (2.15) with $\zeta_i'^n = u_n\Delta Y_{T(n,i)_1}\Delta Y_{T(n,i)_2}1_{\{K(n,i)\geq 2\}}$. So to obtain the result it is enough by Lemma 2.8 to prove that if $\Gamma_t'^n=\sum_{i=1}^{[nt]}\zeta_i'^n$, then $\Gamma_1'^n\to 0$ in law. Then by Lemma 2.5 it is enough to prove that the sequences $(\zeta_i'^n)$ satisfy (2.9) with (2.10). By construction $|\zeta_i'^n|$ is either 0 [with probability $a_n=e^{-\lambda_n}(1+\lambda_n)$] or bigger than $u_n\beta_n^2$, the latter with probability $1-a_n$. Further, $u_n\beta_n^2\to\infty$ in the present case, so (2.9) holds with $\xi_n=\xi_n'=0$ and $\xi_{n,y}''=n(1-a_n)$ for all $n$ large enough. Now $n(1-a_n)\sim\frac{n\lambda_n^2}{2}\to 0$ by (3.13), so indeed we get the result for $j=11$ in Case 1.

Finally, consider $j=12$. Since $|\zeta_i^n(12)|\leq\frac{Cu_n|d_n|}{n}|\Delta Y_{T(n,i)_1}|^2 1_{\{K(n,i)\geq 1\}}$, we readily get

$$E(|\zeta_i^n(12)||\mathcal{F}_{(i-1)/n})\leq\frac{Cu_n|d_n|}{n\theta(\beta_n)}P(K\geq 1)\leq\frac{Cu_n|d_n|}{n^2}.$$

Then (2.7) holds with $\xi_n=Cu_n|d_n|/n$ and we deduce the result for $j=12$ from (3.13).

4.9. *Step* 9: $j=13,14$. First consider for $j=13$. We have

$$E(\zeta_i^n(13)|\mathcal{F}_{(i-1)/n})=\sum_{r=2}^\infty\frac{v_n(r)}{\theta(\beta_n)}(r-1)d_n',$$

$$E(\zeta_i^n(13)^2|\mathcal{F}_{(i-1)/n})\leq C\left(\frac{c_n}{n}+\sum_{r=2}^\infty\frac{v_n(r)}{\theta(\beta_n)}(r-1)^2\right).$$

Therefore, we have (2.8) with $\xi_n=C\lambda_n|d_n'|$ and $\xi_n'=C(c_n+\lambda_n)$ [use (4.12)]. Then (2.11) always holds, and (2.10) holds except in Case 2a: hence, the result for $j=13$.

Finally, we have

$$E(\zeta_i^n(14)|\mathcal{F}_{(i-1)/n})=\frac{d_n}{n}+\frac{d_n'(1-e^{-\lambda_n})}{n\lambda_n}$$

$$= \frac{b'}{n} + \frac{d'_n}{n} \frac{1 - \lambda_n - e^{-\lambda_n}}{\lambda_n},$$

$$E(\zeta_i^n(14)^2 | \mathcal{F}_{(i-1)/n}) \le C\left(\frac{d_n^2}{n^2} + \frac{1}{n}\right).$$

Therefore, we have (2.8) with $\xi_n = C(1 + \lambda_n |d_n|)$ and $\xi'_n = C(1 + d_n^2/n)$. Then (2.10) holds in all cases, and the proof of Lemma 4.1 is now complete.

## 5. Proofs of the theorems.

5.1. *Proof of Theorem* 1.1. By virtue of Corollary 2.3, we can deduce Theorem 1.1 from the tightness of the sequence $(\overline{Y}^n, u_n W^n)$ under Hypothesis (H1-$\alpha$), with our choice of $u_n$ and $\beta_n$. For this, in view of (3.15), it is enough to prove the tightness of the sequence of fourteen-dimensional processes $((\Gamma^n(j)_{1 \le j \le 14})_n$. But this readily follows from Lemmas 2.5 and 4.1, and Theorem 1.1 is proved.

5.2. *Proof of Theorem* 1.2(a). We suppose here that Hypothesis (H2-$\alpha$) holds with $\alpha > 1$ (in particular, we are in Case 1). In view of Lemma 4.1, of (3.15) and of Theorem 2.2, for obtaining Theorem 1.2(a) it suffices to prove that the sequence $(\Gamma^n(14), \Gamma^n(9))$ converges in law to $(Y, W)$, where $W$ is given by (1.7) and $V$ is a Lévy process independent of $Y$ and characterized by (1.8).

Observe that $\zeta_i^n(9)$ satisfies (2.15) with $\zeta_i'^n = u_n \Delta Y_{T(,;i)_1} M_{i/n}^{n,i} 1_{\{K(n,i) \ge 1\}}$. Let $\Gamma_t'^n = \sum_{i=1}^{[nt]} \zeta_i'^n$. In view of Lemma 2.8 [applied with $\eta_i^n = \zeta_i^n(14)$], it is then enough to prove that the pair $(\Gamma^n(14)_1, \Gamma_1'^n)$ converges in law to $(Y_1, V_1)$, where $V_1$ is independent of $Y_1$ and having (1.8) (for $t = 1$). In other words, if we denote by $\phi_n$ the characteristic function of $(\Gamma^n(14)_1, \Gamma_1^n)$, and by $\Phi(u)$ and $\Psi(u)$, respectively, the right-hand sides of (1.2) and (1.8) written for $t = 1$, it suffices to prove that for all $u, v \in \mathbb{R}$, we have

$$\phi_n(u, v) \to \Phi(u)\Psi(v). \tag{5.1}$$

Using (3.17) and (4.25), we get

$$\phi_n(u, v) = e^{iuu_n d_n} \left(e^{-\lambda_n} + \frac{1 - e^{-\lambda_n}}{\theta(\beta_n)} \int_{\{|x| > \beta_n\}} F(dx) e^{iux + z_n(x,v)}\right)^n,$$

where

$$z_n(x, v) = \frac{1}{n} \int_{\{|y| \le \beta_n\}} F(dy)(e^{ivu_n xy} - 1 - iuv_n xy).$$

This is similar to Section 4.7, where $z'_{n,i}$ plays the role of $z_n$ here and, in particular,

$$|z_n(x, v)| \le C \frac{1}{n} u_n^\alpha |v|^\alpha |x|^\alpha. \tag{5.2}$$

We can rewrite $\phi_n$ as

$$\phi_n(u,v) = e^{iuu_n d_n}\Big(1 + \frac{a_n}{n}(A_n(u) + B_n(v) + C_n(u,v))\Big)^n, \tag{5.3}$$

where $a_n = \frac{1-e^{-\lambda_n}}{\lambda_n}$ and

$$A_n(u) = \int_{\{|x|>\beta_n\}} F(dx)(e^{iux} - 1),$$
$$B_n(v) = \int_{\{|x|>\beta_n\}} F(dx)(e^{z_n(x,v)} - 1),$$
$$C_n(u,v) = \int_{\{|x|>\beta_n\}} F(dx)(e^{iux} - 1)(e^{z_n(x,v)} - 1).$$

Combining (5.2) and $|e^{iux}-1| \le |ux|$, and since $x \mapsto |x|^{\alpha+1}$ is $F$-integrable, we first get $|C_n(u,v)| \le C\frac{1}{n}|u||v|^\alpha u_n^\alpha$, hence,

$$C_n(u,v) \to 0. \tag{5.4}$$

Second, $\int_{\{|x|>\beta_n\}} F(dx)(e^{iux} - 1 - iux)$ converges to $\int_{\{|x|>\beta_n\}} F(dx)(e^{iux} - 1 - iux)$, while $d_n = b' - d'_n$, hence,

$$\begin{aligned} A_n(u) + iud_n &\to iub' + \int F(dx)(e^{iux} - 1 - iux) \\ &= iub + \int F(dx)(e^{iux} - 1 - iux\mathbb{1}_{\{|x\le 1\}}). \end{aligned} \tag{5.5}$$

Third, with $K(dx) = \frac{\alpha}{2}((\theta_+^2 + \theta_-^2)\mathbb{1}_{\{x>0\}} + 2\theta_+\theta_-\mathbb{1}_{\{x<0\}})\frac{1}{|x|^{1+\alpha}}\,dx$, we want to prove that

$$B_n(v) \to \int K(dx)(e^{iux} - 1 - iux). \tag{5.6}$$

We have $B_n(v) = B'_n(v) + B''_n(v)$, where

$$B'_n(v) = \int_{\{|x|>\beta_n\}} F(dx) z_n(x,v),$$
$$B''_n(v) = \int_{\{|x|>\beta_n\}} F(dx)(e^{z_n(x,v)} - 1 - z_n(x,v)).$$

Using again (5.2) and the fact that $F$ has compact support and integrates $x \mapsto |x|^{2\alpha}$, we readily see that $|B''_n(v)| \le C|v|^{2\alpha}\frac{u_n^{2\alpha}}{n^2}$, which goes to 0. On the other hand, $B'_n(v) = \int K_n(dx)(e^{ivx} - 1 - ivx)$, where

$$K_n(A) = \frac{1}{n}\int_{\{|x|>\beta_n\}} F(dx) \int_{\{|y|\le\beta_n\}} F(dy)\mathbb{1}_A(u_n xy).$$

Therefore, by Theorem VII.3.4 of Jacod and Shiryaev (2003), (5.6) will follow from the fact that $K_n(h)\to K(h)$ for $h$ equal either to $h_w=\mathbb{1}_{(w,\infty)}$ for $w>0$, or $h'_w=\mathbb{1}_{(-\infty,-w)}$ for $w>0$, or $h'(x)=x^2\mathbb{1}_{\{|x|\le 1\}}$, or $h''(x)=x\mathbb{1}_{\{|x|>1\}}$.

Let us first consider $h=h_w$ for some $w>0$. Since $u_n\beta_n^2\to\infty$ here, for $n$ large enough, $K_n(h_w)$ is the sum of $\gamma_n=\frac{1}{n}\int_{\{x>\beta_n\}}F(dx)F((\frac{w}{u_nx},\beta_n])$ plus another similar term $\gamma'_n$ corresponding to the integrals over the negative half-axis. Further, $\frac{w}{\beta_n u_n x}\to 0$ uniformly in $x>\beta_n$ when $n\to\infty$, therefore, we have

$$\gamma_n\sim\frac{1}{n}\int_{\{x>\beta_n\}}F(dx)\frac{\theta_+u_n^\alpha x^\alpha}{w^\alpha}\sim\frac{\alpha\theta_+^2}{w^\alpha}\frac{u_n^\alpha}{n}\log\frac{1}{\beta_n}\to\frac{\theta_+^2}{2w^\alpha},$$

and, similarly, $\gamma'_n\to\frac{\theta_-^2}{2w^\alpha}$: so $K_n(h_w)\to K(h_w)$. In an analogous fashion we find that $K_n(h'_w)$ converges to $\frac{\theta_+\theta_-}{w^\alpha}=K(h'_w)$. We can also write for $n$ large enough,

$$\begin{aligned}K_n(h')&=\frac{1}{n}\int_{\{|x|>\beta_n\}}F(dx)\int_{\{|y|\le 1/u_n|x|\}}F(dy)u_n^2x^2y^2\\&\sim\frac{1}{n}\int_{\{|x|>\beta_n\}}F(dx)\frac{\alpha\theta}{2-\alpha}u_n^\alpha|x|^\alpha\\&\sim\frac{\alpha^2\theta^2}{2-\alpha}\frac{u_n^\alpha}{n}\log\frac{1}{\beta_n}\to\frac{\alpha\theta^2}{2(2-\alpha)}=K(h').\end{aligned}$$

Finally, we have

$$\begin{aligned}K_n(h'')&=\frac{1}{n}\int_{\{|x|>\beta_n\}}F(dx)\int_{\{1/u_n|x|<|y|\le\beta_n\}}F(dy)u_nxy\\&\sim\frac{1}{n}\int_{\{|x|>\beta_n\}}F(dx)\frac{\alpha\theta'}{\alpha-1}u_n^\alpha|x|^\alpha\\&\sim\frac{\alpha^2\theta\theta'}{\alpha-1}\frac{u_n^\alpha}{n}\log\frac{1}{\beta_n}\to\frac{\alpha\theta\theta'}{2(\alpha-1)}=K(h'').\end{aligned}$$

At this stage we can combine (5.3) with (5.4)–(5.6) to get that $\phi_n(u,v)$ converges to

$$\exp\left(iub+\int F(dx)(e^{iux}-1-iux\mathbb{1}_{\{|x\le 1\}})+\int K(dx)(e^{iux}-1-iux\mathbb{1}_{\{|x\le 1\}})\right),$$

which is $\Phi(u)\Psi(v)$, and we are finished.

5.3. *Proof of Theorem* 1.2(b). We suppose here that Hypothesis (H2-$\alpha$) holds with $\alpha=1$ (in particular, we are in Case 2a). In view of Lemma 4.1, of (3.15) and of Theorem 2.2, for obtaining Theorem 1.2(b), it suffices to prove that the sequence $(\overline{Y}^n,\Gamma^n(10)+\Gamma^n(11))$ converges in probability to $(Y,W)$, where $W_t = -\frac{\theta'^2}{4}\int_0^t g(X_{s-})\,ds$. Since we already know that $\overline{Y}^n \to Y$ (pointwise for the Skorokhod topology) and since $W$ is continuous, it is enough that $\Gamma^n(10) + \Gamma^n(11) \stackrel{P}{\to} W$.

The sum $\zeta_i^n(10)+\zeta_i^n(11)$ satisfies (2.15) with $\zeta_i'^n = \zeta_i'^n(10)+\zeta_i'^n(11)$, where

$$\zeta_i'^n(10) = \frac{u_n d_n^2}{2n^2} + \frac{u_n d_n}{n}\Delta Y_{T(n,i)_1}\mathbb{1}_{\{K(n,i)\geq 1\}},$$
$$\zeta_i'^n(11) = u_n \Delta Y_{T(n,i)_1}\Delta Y_{T(n,i)_2}\mathbb{1}_{\{K(n,i)\geq 2\}}.$$

So by virtue of Lemma 2.8, it is enough to prove that if $\Gamma_t'^n(j) = \sum_{i=1}^{[nt]}\zeta_i'^n(j)$, then

$$\Gamma'^n(10)_1 + \Gamma'^n(11)_1 \stackrel{P}{\to} -\frac{\theta'^4}{4}. \tag{5.7}$$

First, if $a_n = E(\zeta_i'^n(10)|\mathcal{F}_{(i-1)/n})$ and $b_n = E(\zeta_i'^n(10)|\mathcal{F}_{(i-1)/n})$ (which here are nonrandom and independent of $i$), then $a_n$ is given by (4.28) and $b_n$ satisfies (4.29), after replacing the function $g$ by 1. So in view of (3.14) and of $d_n = b' - d_n'$, we get $na_n \to -\frac{\theta'^2}{2}$ and $nb_n \to 0$. Since $E((\Gamma_i'^n(10)_1 - na_n)^2) \leq nb_n$, we get

$$\Gamma'^n(10)_1 \stackrel{P}{\to} -\frac{\theta'^2}{2}. \tag{5.8}$$

Let us use Section 4.8 again: upon replacing once more $g$ by 1, we see that the sequences $(\zeta_i'^n(11))$ satisfy (2.9) with (4.30), hence, $\xi_n' \to 0$ and $\xi_{n,y}'' \to 0$ and $\xi_n \leq C$. So if $\alpha_n = E(\zeta_i'^n(11)\mathbb{1}_{\{|\zeta_i'^n(11)|\leq 1\}}|\mathcal{F}_{(i-1)/n})$ (again nonrandom and independent of $i$), and applying Lemma 2.7, we see that the sequences $(\zeta_i'^n(11) - \alpha_n)$ also satisfy (2.9), with new sequences $\xi_n$, $\xi_n'$ and $\xi_{n,y}''$ all going to 0. Hence, Lemma 2.5 implies $\sum_{i=1}^{[n\cdot]}(\zeta_i'^n(11)-\alpha_n) \stackrel{P}{\to} 0$ and, in view of (5.8), it remains to prove that

$$n\alpha_n \to \frac{\theta'^2}{4}. \tag{5.9}$$

By (3.17) it is clear that

$$\begin{aligned} n\alpha_n &\sim \frac{u_n}{2n}\int_{\{|x|>\beta_n\}} xF(dx)\int_{\{\beta_n<|y|\leq 1/u_n|x|\}} yF(dy) \\ &\sim \frac{1}{2(\log n)^2}\left(d_n'^2 - \int_{\{|x|>\beta_n\}} xF(dx)\int_{\{|y|>1/u_n|x|\}} yF(dy)\right) \end{aligned}$$

for all $n$ large enough, because $F$ has a bounded support and $u_n\beta_n \to 0$. By (3.5), for any $\varepsilon > 0$, there exists $\varepsilon' > 0$ such that $|\frac{d'(\beta)}{\log 1/\beta} - \theta'| \leq \varepsilon$ whenever $\beta \in (0, \varepsilon')$. Then we may write $\int_{\{|x|>\beta_n\}} xF(dx) \int_{\{|y|>1/u_n|x|\}} yF(dy) = x_n + y_n$, where

$$x_n = \int_{\{\beta_n<|x|\leq 1/u_n\varepsilon'\}} xF(dx) \int_{\{|y|>1/u_n|x|\}} yF(dy),$$

$$\begin{aligned} y_n &= \int_{\{|x|>1/u_n\varepsilon'\}} xF(dx) \int_{\{|y|>1/u_n|x|\}} yF(dy) \\ &= \int_{\{|x|>1/u_n\varepsilon'\}} x\, d'\left(\frac{1}{u_n|x|}\right) F(dx). \end{aligned}$$

First, if $|x| \leq \frac{1}{u_n\varepsilon'}$, then $\frac{1}{u_n|x|} \geq \varepsilon'$, so we have

$$|x_n| \leq C_{\varepsilon'} \left| \int_{\{\beta_n<|x|\leq 1/u_n\varepsilon'\}} xF(dx) \right| \leq C_{\varepsilon'}\rho_n \leq C_{\varepsilon'} \log n.$$

Second, if $z_n = \int_{\{|x|>1/u_n\varepsilon'\}} F(dx)|x| \log(u_n|x|)$ and $z'_n = \int_{\{|x|>1/u_n\varepsilon'\}} F(dx)x \times \log(u_n|x|)$, we have $|y_n - \theta' z'_n| \leq \varepsilon z_n$. But (3.6) and (3.5) again imply that

$$\begin{aligned} z_n &= (\log u_n)\left(d_+\left(\frac{1}{u_n\varepsilon'}\right) + d_-\left(\frac{1}{u_n\varepsilon'}\right)\right) + \int_{\{|x|>1/u_n\varepsilon'\}} |x| \log|x| F'(dx) \\ &\sim \theta\left((\log u_n)\log(u_n\varepsilon') - \frac{1}{2}(\log(u_n\varepsilon'))^2\right) \\ &\sim \frac{\theta}{2}(\log n)^2 \end{aligned}$$

and, similarly, $z'_n \sim \frac{\theta'}{2}(\log n)^2$. Hence,

$$n\alpha_n \sim \frac{d'^2_n}{2(\log n)^2} - \frac{x_n + y_n}{2(\log n)^2} \sim \frac{\theta'^2}{2} - \frac{y_n}{2(\log n)^2},$$

and $\frac{y_n}{2(\log n)^2}$ is in between $\frac{\theta' z_n - \varepsilon z'_n}{2(\log n)^2}$ and $\frac{\theta' z_n + \varepsilon z'_n}{2(\log n)^2}$, which respectively converge to $\frac{-\theta'^2 - \theta\varepsilon}{4}$ and to $\frac{-\theta'^2 + \theta\varepsilon}{4}$. Since $\varepsilon > 0$ is arbitrary, we deduce that (5.9) holds, and we are finished.

5.4. *Proof of Theorem* 1.2(c). We suppose here that Hypothesis (H2-$\alpha$) holds for some $\alpha \in (0, 1]$, as well as Hypothesis (H3), and also Hypothesis (H4) if $\alpha < 1$ (so we are in Cases 2b or 3b). In view of Lemma 4.1, of (3.15) and of Theorem 2.2, it suffices to prove that the sequence $(\Gamma^n(14), \Gamma^n(11))$ converges in law to $(Y, W)$, where $W$ is given by (1.7) and $V$ is a Lévy process independent of $Y$ and characterized by (1.10).

1. Note that if

$$\Gamma'^n_t = \sum_{i=1}^{[nt]} \zeta'^n_i \qquad \text{where } \zeta^n_i = u_n \Delta Y_{T(n,i)_1} \Delta Y_{T(n,i)_2} \mathbb{1}_{\{K(n,i)\geq 2\}}, \tag{5.10}$$

$$(5.11)\quad Z_t^n = \sum_{i=1}^{[nt]} \eta_i^n \qquad \text{where } \eta_i^n = \Delta Y_{T(n,i)_1} \mathbb{1}_{\{K(n,i)\geq 1\}},$$

then first $\zeta_i^n(11)$ satisfies (2.15) with $\zeta_i'^n$ as above, and $\Gamma^n(14) - Z^n$ converges locally uniformly in time to $bt$ (because $d_n = b$ here). So in view of Lemma 2.8, it is enough to prove that $(Z_1^n, \Gamma_1'^n)$ converges in law to $(Z_1, V_1)$, where $Z_t = Y_t - bt$.

2. For each $n$ the variables $(\eta_i^n, \zeta_i'^n)_{i\geq 1}$ are i.i.d. centered, and we denote by $K_n$ their law, which is given by

$$K_n(A) = \frac{e^{-\lambda_n}\lambda_n}{\theta(\beta_n)} \int_{\{|x|>\beta_n\}} F(dx) \mathbb{1}_A(x,0) + \frac{1-e^{-\lambda_n}(1+\lambda_n)}{\theta(\beta_n)^2} \int_{\{|x|>\beta_n\}} F(dx) \int_{\{|y|>\beta_n\}} F(dy) \mathbb{1}_A(x, u_n xy).$$

Since the two processes $Z$ and $V$ are independent Lévy processes, they have no common jumps and, further, the jumps of $Z$ are the same as those of $Y$: so the Lévy measure $K$ of the pair $(Z,V)$ is

$$K(dx,dy) = F(dx)\varepsilon_0(dy) + \frac{\theta^2\alpha}{4}\varepsilon_0(dx)\frac{1}{|y|^{1+\alpha}}\,dy.$$

Observe that $(Z,V)$ has no drift and no continuous martingale part.

By virtue of Theorem VII.3.4 of Jacod and Shiryaev (2003), for the convergence in law of $(Z_1^n, \Gamma_1'^n)$ to $(Z_1, V_1)$, it is enough to prove that

$$(5.12)\qquad nK_n(h) \to K(h)$$

for all $h$ which are continuous bounded and vanish on a neighborhood of 0, and also for $h = h_u, h_u', h_u''$, where $h_u(x,y) = x^2 \mathbb{1}_{\{x^2+y^2\leq u\}}$ and $h_u'(x,y) = y^2 \mathbb{1}_{\{x^2+y^2\leq u\}}$ and $h_u''(x,y) = xy\mathbb{1}_{\{x^2+y^2\leq u\}}$, for almost all $u > 0$ (for Lebesgue measure). Since both $K$ and $K_n$ are invariant under the maps $(x,y) \mapsto (-x,y)$ and $(x,y) \mapsto (-x,y)$, it is even enough to prove (5.12) for $h_u$, $h_u'$, $h_u''$, and also for $h_{u,v}(x,y) = \mathbb{1}_{\{|x|\geq u, |y|\geq v\}}$ for all $u, v \in \mathbb{R}_+$ such that $(u,v) \neq (0,0)$.

We begin with (5.12) for $h_{u,v}$. For $u > 0$, we have as soon as $\beta_n < u$,

$$nK_n(h_{u,0}) = \frac{ne^{-\lambda_n}\lambda_n}{\theta(\beta_n)}\theta(u-) + n\frac{1-e^{-\lambda_n}(1+\lambda_n)}{\theta(\beta_n)^2} \int_{\{|x|\geq u\}} F(dx)\theta(\beta_n),$$

where $\theta(u-)$ denotes the left limit at point $u$ of the decreasing and right-continuous function $\theta(\cdot)$. The last term above is smaller than $C_u\theta(\beta_n)/n$, which goes to 0 as $n \to \infty$, and the first term converges to $\theta(u-) = K(h_{u,0})$, so (5.12) holds for $h_{u,0}$.

Next, if $u>0$ and $v>0$, we have as soon as $\beta_n<u$ (recall that $u_n=1/\beta_n$ here),

$$\begin{aligned} nK_n(h_{u,v}) &= n\frac{1-e^{-\lambda_n}(1+\lambda_n)}{\theta(\beta_n)^2}\int_{\{|x|\geq u\}}F(dx)\int_{\{|y|>\beta_n,|y|\geq\beta_n v/|x|\}}F(dy)\\ &\leq C_u\frac{\theta(\beta_n)}{n}, \end{aligned}$$

which again goes to 0: so (5.12) holds for $h_{u,v}$. Finally, if $v>0$, we have as soon as $\beta_n<v$,

$$\begin{aligned} nK_n(h_{0,v}) &= n\frac{1-e^{-\lambda_n}(1+\lambda_n)}{\theta(\beta_n)^2}\int_{\{|x|>\beta_n\}}F(dx)\int_{\{|y|>\beta_n,|y|\geq\beta_n v/|x|\}}F(dy)\\ &\sim\frac{1}{2n}\int_{\{|x|>\beta_n\}}F(dx)\int_{\{|y|>\beta_n,|y|\geq\beta_n v/|x|\}}F(dy). \end{aligned}\tag{5.13}$$

Let $\varepsilon>0$. By Hypothesis (H2-$\alpha$) there exists $\varepsilon'\in(0,v)$ such that $|\beta^\alpha\theta(\beta)-\theta|\leq\varepsilon$ for all $\beta\in(0,2\varepsilon']$. By (5.13) we see that $nK_n(h_{0,v})\sim x_n+y_n+z_n$, where

$$\begin{aligned} x_n &= \frac{1}{2n}\int_{\{\beta_n<|x|\leq v\beta_n/\varepsilon'\}}F(dx)\int_{\{|y|\geq\beta_n v/|x|\}}F(dy),\\ y_n &= \frac{1}{2n}\int_{\{v\beta_n/\varepsilon'<|x|\leq v\}}F(dx)\int_{\{|y|\geq\beta_n v/|x|\}}F(dy),\\ z_n &= \frac{1}{2n}\int_{\{|x|>v\}}F(dx)\int_{\{|y|>\beta_n\}}F(dy). \end{aligned}$$

Using Hypothesis (H2-$\alpha$) again and (3.2), we get

$$z_n\leq C_v\frac{\theta(\beta_n)}{n},\tag{5.14}$$

$$\begin{aligned} x_n &\leq \frac{C}{n\beta_n^\alpha v^\alpha}\int_{\{\beta_n<|x|\leq v\beta_n/\varepsilon'\}}|x|^\alpha F(dx)\\ &= \frac{C}{n\beta_n^\alpha v^\alpha}\int_0^{v\beta_n/\varepsilon'}y^{\alpha-1}\theta(y\vee\beta_n)(dy)\\ &\leq \frac{C_v}{n\beta_n^\alpha}\left(1+\log\frac{v}{\varepsilon}\right)\leq\frac{C_{v,\varepsilon}}{n\beta_n^\alpha}. \end{aligned}\tag{5.15}$$

Finally, if $y'_n=\frac{\theta}{2n\beta_n^\alpha v^\alpha}\int_{\{v\beta_n/\varepsilon<|x|\leq v\}}|x|^\alpha F(dx)$, we have $y'_n(1-\varepsilon')\leq y_n^\varepsilon\leq y'_n(1+\varepsilon')$. But $y'_n=\frac{\theta}{2n\beta_n^\alpha v^\alpha}(\rho(\beta_n v/\varepsilon)-\rho(v))$, which is equivalent to $\frac{\alpha\theta^2}{2n\beta_n^\alpha v^\alpha}\log\frac{1}{\beta_n}$ by (3.5). Putting this together with (5.14) and (5.15), and since $\varepsilon$ is arbitrarily small, we deduce that $nK_n(h_{0,v})$ is equivalent to $\frac{\alpha\theta^2}{2n\beta_n^\alpha v^\alpha}\log\frac{1}{2n\beta_n}$. Then, clearly, $nK_n(h_{0,v})$ converges to $\frac{\theta^2}{2v^\alpha}$, which equals $K(h_{0,v})$.

It remains to prove (5.12) for $h_u$, $h'_u$ and $h''_u$. First, because of Hypothesis (H3) we have $K_n(h''_u)=0$ and $K(h''_u)=0$. Next,

$$\begin{aligned}
nK_n(h_u) = {} & \frac{ne^{-\lambda_n}\lambda_n}{\theta(\beta_n)}\int_{\{\beta_n<|x|\le\sqrt{u}\}} x^2F(dx)\\
& + n\frac{1-e^{-\lambda_n}(1+\lambda_n)}{\theta(\beta_n)^2}\\
& \times\int_{\{\beta_n<|x|<\sqrt{u}\}} x^2F(dx)\int_{\{\beta_n<|y|\le\beta_n\sqrt{u/x^2-1}\}} F(dy).
\end{aligned}$$

The last term above is smaller than $C\theta(\beta_n)/n$, which goes to 0, while the first term goes to $\int_{\{|x|\le\sqrt{u}\}} x^2F(dx)$, which equals $K(h_u)$, hence, (5.12) holds for $h=h_u$. Finally, we have as soon as $\beta_n<\sqrt{u}$,

$$\begin{aligned}
nK_n(h'_u) &= u_n^2 n\frac{1-e^{-\lambda_n}(1+\lambda_n)}{\theta(\beta_n)^2}\\
&\quad\times\int_{\{\beta_n<|x|<\sqrt{u/2}\}} x^2F(dx)\int_{\{\beta_n<|y|\le\beta_n\sqrt{u/x^2-1}\}} y^2F(dy)\\
&\sim \frac{u_n^2}{2n}\int_{\{\beta_n<|x|<\sqrt{u/2}\}} x^2c\left(\beta_n\sqrt{u/x^2-1}\right)F(dx)\\
&\quad-\frac{u_n^2c_n}{2n}\int_{\{\beta_n<|x|<\sqrt{u/2}\}} x^2F(dx).
\end{aligned}$$

The last term above is smaller that $Cu_n^2c_n/n$, which goes to 0. In view of (3.3), the first term is equivalent to

$$\begin{aligned}
&\frac{u_n^2\beta_n^{2-\alpha}\alpha\theta}{2(2-\alpha)n}\int_{\{\beta_n<|x|<\sqrt{u/2}\}} x^2\left(\frac{u}{x^2}-1\right)^{(2-\alpha)/2}F(dx)\\
&\qquad\sim\frac{u_n^\alpha\alpha\theta u^{(2-\alpha)/2}}{2(2-\alpha)n}\int_{\{\beta_n<|x|<\sqrt{u/2}\}}|x|^\alpha F(dx)\\
&\qquad\sim\frac{u_n^\alpha\alpha^2\theta^2u^{(2-\alpha)/2}}{2(2-\alpha)n}\log\left(\frac{1}{\beta_n}\right),
\end{aligned}$$

which converges to $\frac{\alpha\theta^2u^{(2-\alpha)/2}}{2(2-\alpha)}$, which in turn equals $K(h'_u)$: we are thus finished.

REMARK 5.1. When Hypothesis (H2-$\alpha$) holds for some $\alpha\in(1,2)$ and also Hypothesis (H3) holds, one could prove part (a) of Theorem 1.2 by the same method as above for (c): we have, of course, the same $u_n=(n/\log n)^{1/\alpha}$, but instead of $\beta_n=\log n/n^{1/2\alpha}$, one could take $\beta_n=(\log n/n)^{1/\alpha}$ as in Cases 2b and 3b. Then in Lemma 4.1 one obtains that $\Gamma^n(j)$ goes to 0 for $j=9$, but not for $j=11$.

5.5. *Proof of Theorem* 1.2(d). Now we assume Hypothesis (H1-$\alpha$) for some $\alpha<1$ (i.e., Case 3a). In view of Lemma 4.1, of (3.15) and of Theorem 2.2, for obtaining Theorem 1.2(d) it suffices to prove that the sequence $(\Gamma^n(14),\Gamma^n(10)+\Gamma^n(12))$ converge in law to $(Y,W)$, where $W$ is given by (1.11), which can also be written as

$$\begin{aligned}W_t = d\sum_{n:\,R_n\le t}&([f(X_{R_n-}+\Delta Y_{R_n}f(X_{R_n-}))-f(X_{R_n-})]\xi_n\\&\qquad + f(X_{R_n-})f'(X_{R_n-})\Delta Y_{R_n}(1-\xi_n))\\&+\frac{d^2}{2}\int_0^t f(X_{s-})f'(X_{s-})\,ds.\end{aligned}\tag{5.16}$$

1. We have $\zeta_i^n(10)=\zeta_i'^n(10)+\zeta_i''^n(10)$, where $\zeta_i''^n(10)=u_n g(X_{(i-1)/n})\frac{d_n^2}{2n^2}$ and

$$\zeta_i'^n(10)=u_n g(X_{(i-1)/n})\frac{d_n}{n}\Delta Y_{T(n,i)_1}1_{\{K(n,i)\ge1\}}.$$

Set $\Gamma_t'^n=\sum_{i=1}^{[nt]}\zeta_i'^n(10)$ and $\Gamma_t''^n=\sum_{i=1}^{[nt]}\zeta_i''^n(10)$. Observe that $\frac{u_n d_n^2}{n}\to d^2$, hence, a simple Riemann approximation of the Lebesgue integral shows the following convergence, locally uniform in $t$:

$$\Gamma'^n(10)_t\stackrel{P}{\to}\frac{d^2}{2}\int_0^t g(X_{s-})\,ds.\tag{5.17}$$

Therefore, it is enough to prove that the pair $(\Gamma^n(14),\Gamma''^n+\Gamma^n(12))$ converges in law to the pair $(Y,\Gamma')$, where $\Gamma'$ is the first term in the right-hand side of (1.11).

Set $a_n=u_nd_n/n$. We can write $\Gamma''^n+\Gamma^n(12)=a_n\Gamma^n$, where $\Gamma_t^n=\sum_{i=1}^{[nt]}\zeta_i^n$ and

$$\begin{aligned}\zeta_i^n=(&G(X_{(i-1)/n},\Delta Y_{T(n,i)_1})(i-nT(n,i)_1)\\&+g(X_{(i-1)/n})\Delta Y_{T(n,i)_1}(nT(n,i)_1-i+1))1_{\{K(n,i)\ge1\}}.\end{aligned}\tag{5.18}$$

We also write $V'=d\Gamma$, where [cf. (1.11)]

$$\begin{aligned}\Gamma_t=\sum_{n:\,R_n\le t}&([f(X_{R_n-}+\Delta Y_{R_n}f(X_{R_n-}))-f(X_{R_n-})]\xi_n\\&\qquad + g(X_{R_n-})\Delta Y_{R_n}(1-\xi_n)).\end{aligned}\tag{5.19}$$

Since $a_n \to d$, it remains to prove that $(\Gamma^n(14), \Gamma^n)$ converges to $(Y, \Gamma)$ in law. Observe also that $d_n \to d$, so exactly as in Section 6, and if we set (5.11), it is enough to prove that $(Z^n, \Gamma^n)$ converges in law towards $(Z, \Gamma)$, where $Z_t = Y_t - td$.

2. The presence of an apriori arbitrary function $G$ in (5.18) makes things a bit difficult, and in the absence of a general theory to handle this case, we use a trick, pretending first that $X_{(i-1)/n}$ does not show up in (5.18). That is, with an arbitrary measurable bounded function $l$ on $\mathbb{R} \times [0,1]$, with support in $[-p,p] \times [0,1]$, which satisfies $|l(x,u)| \le C|x|$, we set

$$\Gamma^n(l)_t = \sum_{i=1}^{[nt]} \zeta_i^n(l) \tag{5.20}$$
$$\text{where } \zeta_i^n(l) = l(\Delta Y_{T(n,i)_1}, nT(n,i)_1 - i + 1)\mathbb{1}_{\{R(n,i) \ge 1\}}.$$

We will study the convergence of the pair $(Z^n, \Gamma^n(l))$. In view of (3.17), the law $K_n$ of the pair $(\eta_i^n, \zeta_i^n(l))$ is independent of $i$ and given by

$$\begin{aligned} K_n(h) &= e^{-\lambda_n} h(0,0) \\ &\quad + \frac{1}{n} \int_0^1 e^{-\lambda_n u}\, du \int_{\{|x| > \beta_n\}} F(dx) \int_0^1 h(x, l(x,u))\, du. \end{aligned} \tag{5.21}$$

Letting $K$ be the measure on $\mathbb{R}^2$ defined by

$$K(h) = \int F(dx) \int_0^1 h(x, l(x,u))\, du, \tag{5.22}$$

we want to prove that $nK_n(h) \to K(h)$ [i.e., (5.12)] holds for suitable functions $h$.

Suppose first that $h$ is continuous and bounded and vanishes on a neighborhood of 0. Observe that since $|l(x,u)| \le C|x|$, we have $h(x, l(x,u)) = 0$ if $|x| \le \varepsilon$ for some $\varepsilon > 0$. In (5.21) the contribution of the first term to the right to $nK_n(h)$ is 0; as soon as $\beta_n < \varepsilon$, the contribution of the second term is

$$\int_0^1 e^{-\lambda_n u}\, du \int_{\{|x| > \varepsilon\}} F(dx) \int_0^1 h(x, l(x,u)),$$

which obviously converge to $K(h)$ because $\lambda_n \to 0$: so we have (5.12).

Now we take the function $h(x,y) = x\mathbb{1}_{\{|x| \le v\}}$ for any given $v > 0$. In (5.21) the contribution to $nK_n(h)$ of the first term to the right is 0, and the contribution of the second term is

$$\int_0^1 e^{-\lambda_n u}\, du \int_{\{\beta_n < |x| \le v\}} xF(dx),$$

which goes to $K(h)$ again. The same argument works as well for the functions $h(x,y)=y\mathbb{1}_{\{|x|\le v\}}$, $h(x,y)=x^2\mathbb{1}_{\{|x|\le v\}}$, $h(x,y)=xy\mathbb{1}_{\{|x|\le v\}}$ and $h(x,y)=y^2\mathbb{1}_{\{|x|\le v\}}$ [here we do not truncate in $y$, since the argument $y$ is replaced by the bounded term $l(x,u)$ in (5.21)].

All this shows that, by using Theorem VII.3.4 of Jacod and Shiryaev (2003), the pair $(Z^n,\Gamma^n(l))$ converges in law to a Lévy process $(Z,\Gamma(l))$ with no drift, no continuous part and Lévy measure $K$, and also has (UT) by Lemma 2.1.

We can, of course, realize this pair $(Z,\Gamma(l))$ as follows: we first take $Z_t=Y_t-dt$ (so this is in accordance with our previous notation), and we label the jump times of $Y$ as $R_1,R_2,\dots$ [as in (5.19)]. Then we define, possibly on an extension of the space, a sequence $(\xi_n)$ of i.i.d. variables, independent of $Y$, and uniform over $[0,1]$. Then we set

$$\Gamma(l)_t=\sum_{n:\,R_n\le t} l(\Delta Y_{R_n},\xi_n). \tag{5.23}$$

Observe that $\Gamma(l)$, as well as $Y$ and $Z$, has finite variation over finite intervals.

3. We will apply the preceding results to the functions $l=l_y$ defined by $l_y(x,u)=G(y,x)(1-u)+g(y)xu$: we call $\Gamma^n(y)$ and $\Gamma(y)$ the processes $\Gamma^n(l_y)$ and $\Gamma(l_y)$, and write also $\zeta_i^n(y)=\zeta_i^n(l_y)$. If we pick finitely many arbitrary points $y_j$, not only do we have the convergence in law of $(Z^n,\Gamma^n(y_j))$ to $(Z,\Gamma(y_j))$ for any given $j$, but one could prove in a similar way that we have the convergence of $(Z^n,\Gamma^n(y_1),\dots,\Gamma^n(y_k))$ to $(Z,\Gamma(y_1),\dots,\Gamma(y_k))$, with the same sequence $\xi_n$ in the definition of all $\Gamma(y_j)$'s: one just has to put any finitely many functions $l$'s in (5.21) and (5.22) to see that everything works out the same way. This, of course, gives the convergence for an infinite sequence of $y_k$'s.

So we pick a dense sequence $y_k$. By the Skorokhod representation theorem we can find another probability space on which new processes still called $(Y,Z,\Gamma(y_1),\dots,\Gamma(y_k),\dots)$ and $(Y,Z^n,\Gamma^n(y_1),\dots,\Gamma^n(y_k)\dots)$ are defined and have the same law as the original ones and, further, $(Z^n,\Gamma^n(y_1),\dots,\Gamma^n(y_k)\dots)$ converges pointwise for the Skorokhod topology on $\mathbb{D}(\mathbb{R}_+,\mathbb{R}^{\mathbb{N}})$ to $(Z,\Gamma(y_1),\dots,\Gamma(y_k),\dots)$. On the new space we still have the representation (5.23) for all $\Gamma(y_k)$ with the same sequence $(\xi_n)$. Furthermore, on the new space we can solve our equation (1.1), having a solution $X$, and redefine $\zeta_i^n$ by (5.18) and $\Gamma_t^n=\sum_{i=1}^{[nt]}\zeta_i^n$: here $Y$ has locally finite variation, so the filtrations play no role at all. So we have $\zeta_i^n=l_{X_{(i-1)/n}}(\Delta Y_{T(n,i)_1},nS_1)\mathbb{1}_{\{K(n,i)\le 1\}}$.

Now the functions $y\mapsto l_y(x,u)$ are continuous, and even much more. Namely, we have

$$w_K(\varepsilon):=\sup\biggl(\frac{|l_y(x,u)-l_{y'}(x,u)|}{|y|}:y,y'\in K,|y-y'|\le\varepsilon,$$

$$x \in [-p,p], u \in [O,1]\Big) \to 0$$

as $\varepsilon \to 0$ for any compact set $K$. Let us first pick a point $\omega$ in the event space and a time $t > 0$, then a compact $K$ containing the path of $X$ over $[0,t]$, and an $\varepsilon > 0$. We can find a finite subdivision $t_0 = 0 < t_1 < \cdots < t_k = t$ and points $y_j \in K$, such that $|X_t - y_j| \le \varepsilon$ for all $t \in [t_{j-1}, t_j)$. Set $J(n,j) = \{i : \frac{i-1}{n} \in [t_{j-1}, t_j)\}$. If $i \in J(n,j)$, we have $|\zeta_i^n - \zeta_i^n(y_j)| \le w_K(\varepsilon)|\Delta Y_{S_1}|$. It follows that for all $s \le t$,

$$\left|\Gamma_s^n - \sum_{j=1}^k (\Gamma^n(y_j)_{s\wedge t_j} - \Gamma^n(y_j)_{s\wedge t_{j-1}})\right| \le w_K(\varepsilon) \sum_{r\le t} |\Delta Y_r|. \tag{5.24}$$

Similarly, if we set

$$\Gamma_t = \sum_{n:\, R_n \le t} l_{X_{R_n-}}(\Delta Y_{R_n}, U_n) \tag{5.25}$$

[i.e., $\Gamma$ is defined by (5.19), on our new space], we get for all $s \le t$,

$$\left|\Gamma_s - \sum_{j=1}^k (\Gamma'(y_j)_{s\wedge t_j} - \Gamma(y_j)_{s\wedge t_{j-1}})\right| \le w_K(\varepsilon) \sum_{r\le t} |\Delta Y_r|. \tag{5.26}$$

At this point, it suffices to use that $(Z^n, \Gamma^n(y_k), \dots, \Gamma^n(y_k)) \to (Z, \Gamma(y_k), \dots, \Gamma(y_k))$ for the Skorokhod topology to obtain that the upper limit of the Skorokhod distance between $(Z^n, \Gamma^n)$ and $(Z, \Gamma)$ over $[0,t]$ is smaller or equal to $w_K(\varepsilon) \sum_{s\le t} |\Delta Y_s| < \infty$. Since $w_K(\varepsilon) \to 0$ as $\varepsilon \to 0$ and since $\sum_{s\le t} |\Delta Y_s| < \infty$, we conclude that for our particular point $\omega$ we have $(Z^n, \Gamma^n) \to (Z, \Gamma)$.

This works for all points $\omega$. Going back to the original space, this clearly implies that indeed $(Z^n, \Gamma^n)$ converges in law to $(Z, \Gamma)$, and we are finished.

Laboratoire de Probabilités
et Modèles Aléatoires
CNRS UMR 7599
Université Pierre et Marie Curie
Tour 56
4 Place Jussieu
75252 Paris Cedex 05
France
e-mail: jj@ccr.jussieu.fr